\documentclass[11pt,a4paper,leqno]{amsart}
\usepackage{graphicx} 
\usepackage{amssymb}
\usepackage{amsmath}
\usepackage{amsthm}
\usepackage{cancel}
\usepackage{amsfonts,amssymb}
\usepackage[abbrev,non-sorted-cites, alphabetic]{amsrefs}
\usepackage{latexsym}
\usepackage{color}
\usepackage{appendix}
\usepackage{hyperref}
\hypersetup{
    colorlinks=true,
    allcolors=black,
}

\usepackage{xr}
\oddsidemargin=.in
\evensidemargin=.in

\newtheorem{theorem}{Theorem}[section]
\newtheorem{lem}[theorem]{Lemma}
\theoremstyle{plain}

\newtheorem{assumption}[theorem]{Assumption}

\newtheorem{claim}{Claim}

\newtheorem{corollary}[theorem]{Corollary}

\newtheorem{definition}[theorem]{Definition}

\newtheorem{lemma}[theorem]{Lemma}

\newtheorem{proposition}[theorem]{Proposition}
\newtheorem{remark}{Remark}

\newenvironment{claimproof}{\par \noindent \textit{Proof of Claim.}}{\hfill $\diamond$ \par}

\newcommand{\N}{\mathbb{N}}     
\newcommand{\R}{\mathbb{R}}     
\newcommand{\C}{\mathbb{C}}     

\newcommand{\Cone}{\mathcal{C}} 

\newcommand{\loc}{\mathrm{loc}}       

\title[LMCF Singularities with Bounded $H$]{Existence of Lagrangian Mean Curvature Flow Singularities with Bounded Mean Curvature}

\author{Maxwell Stolarski}
\address{Maxwell Stolarski: School of Mathematics and Statistics, University College Dublin, Science Centre North, Belfield, Dublin 4, Ireland}
\email{maxwell.stolarski@ucd.ie}

\author{Wei-Bo Su}
\address{Wei-Bo Su: Department of Mathematics, National Central University, No. 300, Zhongda Rd., Zhongli District, Taoyuan City 320317, Taiwan}
\email{weibosu@math.ncu.edu.tw}

\date{\today}

\begin{document}
    
\begin{abstract}
In every dimension $n \ge 3$, we show there exist Lagrangian mean curvature flows that develop neck pinch singularities such that the mean curvature remains uniformly bounded.
This follows from a related result that, under a structural assumption on the singularity formation, the mean curvature at small scales near the singularity can be bounded by the Lagrangian angle at outer scales with suitable weights.
Applying this result to a subset of the finite-time Lagrangian mean curvature flow singularities constructed in \cite{SS26II} and the infinite-time Lagrangian mean curvature flow singularities constructed in \cite{STW} then gives examples of both finite-time and infinite-time singularity formation with uniformly bounded mean curvature.
\end{abstract}

\maketitle


\section{Introduction}

Mean curvature flow $\{F_t:L\to L(t)\}_{t\in I}$ evolves a submanifold by its mean curvature vector $\mathbf{H}_{L(t)}$. As the negative $L^2$-gradient flow of volume, it provides a natural framework for studying minimal submanifolds. However, the flow may develop singularities at which the second fundamental form $\mathbf{A}_{L(t)}$ becomes unbounded \cite{Huisken84}:
\[
    \limsup_{t\to T}\sup_{L(t)}|\mathbf{A}_{L(t)}|=\infty,
\]
where the singular time $T$ may be finite or infinite. A finite-time singularity is called Type I if
\[
    \sup_{L(t)}|\mathbf{A}_{L(t)}|
    \leq \frac{C}{\sqrt{T-t}}
\]
for some constant $C>0$ and all $t<T$ sufficiently close to $T$; otherwise, it is called Type II.

Type I singularities are modeled on smooth, non-flat self-shrinkers under parabolic rescaling, and the mean curvature necessarily becomes unbounded as well \cite{LeSesum10}. By contrast, certain Type II singularities are modeled on non-flat minimal submanifolds at scales smaller than the parabolic scale $\sqrt{T-t}$. This raises the possibility that the second fundamental form may blow up while the mean curvature remains uniformly bounded.

For mean curvature flow of closed embedded surfaces in $\mathbb{R}^{3}$, Li--Wang \cite{LiWang19} proved the mean curvature must blow up at the first singular time, regardless of whether the singularity is Type I or Type II. Le-{\v S}e{\v s}um \cite{LeSesum10} have an alternative proof of this result assuming the multiplicity-one conjecture, which was recently proven by Bamler-Kleiner \cite{BamlerKleiner23}. In higher dimensions, Vel\'azquez \cite{Velazquez94} constructed a family of complete, $O(n)\times O(n)$-invariant hypersurfaces in $\mathbb{R}^{2n}$, $n\geq 4$, whose mean curvature flows develop finite-time Type II singularities with precisely quantified blow-up rates. At the curvature scale, these singularities are modeled on smooth minimal hypersurfaces asymptotic to the Simons cone. The blow-up rates are governed by eigenvalues of the linearized rescaled mean curvature flow operator at this cone. Guo--\v{S}e\v{s}um \cite{GuoSesum18} showed that, for $n\geq 5$, solutions corresponding to the lowest admissible mode exhibit mean curvature blow-up, albeit at a slower rate than the second fundamental form. In contrast, the first-named author proved solutions corresponding to higher even modes have uniformly bounded mean curvature, despite the faster blow-up of their second fundamental form \cite{Stolarski23}.

Analogous questions apply to the Ricci flow.
At finite-time Ricci flow singularities, the Riemann curvature $Rm$ and Ricci curvature $Rc$ always blow up \cites{Hamilton95, Sesum05}.
The behavior of scalar curvature $S$, however, is more subtle.
The first-named author constructed Ricci flow solutions that form finite-time singularities modeled on Ricci-flat cones \cite{Stolarski24}, and
OpenAI very recently indicated these solutions have uniformly bounded scalar curvature $S$ in high dimensions \cite{OpenAI_Bdd_S}. 
In low dimensions $n \le 3$, the Hamilton-Ivey pinching estimates \cites{Hamilton95, Ivey93} imply $S$ always blows up at finite-time Ricci flow singularities.
OpenAI have also claimed $S$ always blows up at finite-time Ricci flow singularities in dimension $n =4$ \cite{OpenAI_S_4D}.

In this paper, we study singularities of Lagrangian mean curvature flow (LMCF), the mean curvature flow of Lagrangian submanifolds in Calabi--Yau manifolds. The conjectural programme initiated by Thomas--Yau \cite{ThomasYau02} and refined by Joyce \cite{Joyce15} proposes using LMCF, with surgeries when necessary, to construct special Lagrangian submanifolds, which minimize volume in their homology classes. Understanding singularity formation is therefore essential to this programme.

A fundamental example of finite-time singularity formation is the Lawlor neck pinch, whose tangent flow is the static special Lagrangian union of two transverse planes and whose curvature-scale blow-up model is a Lawlor neck, a special Lagrangian desingularization of this union. Neves \cite{Neves07} established finite-time singularity formation in an equivariant setting, while Wood \cite{Wood24} identified the corresponding tangent flows and Lawlor neck blow-up models. A major advance beyond the equivariant setting was made by Lotay--Schulze--Sz\'ekelyhidi \cite{LSS22}. For zero-Maslov, rational LMCF in compact Calabi--Yau surfaces, they proved uniqueness of the two-plane tangent flow and established a detailed geometric description of the collapsing Lawlor neck, without symmetry assumptions. Furthermore, when all singularities at the first singular time are of this form, they showed that the flow can be continued past the singular time as a smooth immersed LMCF. In arbitrary dimension, Li--Sz\'ekelyhidi \cite{LS24} subsequently proved, under suitable hypotheses, that a singularity admitting such a two-plane tangent flow has a Type II blow-up sequence converging to a Lawlor neck. Székelyhidi \cite{Szekelyhidi26} introduced a notion of nondegenerate neck pinch for LMCF of surfaces, proved that such singularities are stable under small perturbations, and constructed examples in compact Calabi--Yau surfaces. He also showed that neck pinches can be perturbed near the singular time to become nondegenerate. The authors constructed finite-time Type II singularities, including Lawlor neck pinches, with quantitatively precise curvature blow-up rates \cite{SS26II}. Infinite-time Lawlor neck pinches with precise polynomial curvature blow-up rates were constructed by the second-named author with Tsai and Wood \cite{STW}.

\subsection{Main Result}
Our main result is a local estimate for the mean curvature near
Lawlor neck pinches in complex dimensions $n\geq 3$, under the
following two-scale structural assumption.

\begin{assumption}\label{blow-up assumptions}
    Let $0<T\leq\infty$, and let
    $\{F_t:L\to L(t)\}_{t\in[0,T)}$ be a smooth, properly immersed,
    zero-Maslov LMCF in $\C^n$, with $n\geq 3$.
    Suppose that there exist an \emph{outer scale}
    $R:[0,T)\to(0,r_0]$ and an \emph{inner scale}
    $\lambda:[0,T)\to(0,\infty)$ satisfying the following conditions.
    \begin{enumerate}
        \item\label{scale assumptions}
        Both $R$ and $\lambda$ are $C^1$, $R$ is non-increasing,
        and the following scale comparison estimates hold:
        \begin{enumerate}
            \item\label{scale assumption 1}
            the inner and outer scales are asymptotically separated
            \begin{equation}\label{eqn scale assumption 1}
                \lim_{t\to T}\frac{\lambda(t)}{R(t)}=0,
            \end{equation}
            and
            \item\label{scale assumption 2}
            there exists $D_0>0$ such that
            \begin{equation}\label{eqn scale assumption 2}
                R(t)^2\frac{|\lambda'(t)|}{\lambda(t)}
                \leq D_0
            \end{equation}
            for all $t\in[0,T)$.
        \end{enumerate}

        \item\label{assumption, blow-up to Lawlor neck}
        For every sequence $t_i\to T$, the parabolically rescaled flows
        \[
            L_i(\tau)
            :=\lambda(t_i)^{-1}
              L\bigl(t_i+\tau\lambda(t_i)^2\bigr)
        \]
        admit a subsequence converging in
        $C^\infty_{\mathrm{loc}}
        \bigl(\C^n\times(-\infty,0]\bigr)$
        to the static flow generated by a multiplicity one Lawlor neck $\overline L$.

        \item\label{assumption, blow-up to cone}
        For every pair of sequences $t_i\to T$ and $\mu_i\to 0$
        satisfying
        \[
            \frac{\lambda(t_i)}{\mu_i}\to 0
            \qquad \text{and} \qquad 
            \frac{\mu_i}{R(t_i)}\to 0,
        \]
        the parabolically rescaled flows
        \[
            L_i(\tau)
            :=\mu_i^{-1}L\bigl(t_i+\tau\mu_i^2\bigr)
        \]
        admit a subsequence converging in
        $C^\infty_{\mathrm{loc}}
        \bigl((\C^n\setminus\{\mathbf 0\})\times(-\infty,0]\bigr)$
        to the static flow generated by a multiplicity one special Lagrangian cone
        $\mathcal C=\mathcal P_0\cup\mathcal P_1$,
        where $\mathcal P_0$ and $\mathcal P_1$ are two
        transversely intersecting planes.
    \end{enumerate}
\end{assumption}

Bounding the mean curvature $\mathbf H$ under these assumptions relies on the evolution equation
$\partial_t\theta_{L(t)}=\Delta_{L(t)}\theta_{L(t)}$
for the Lagrangian angle under LMCF and the identity
$\mathbf H_{L(t)}=J\nabla\theta_{L(t)}$.
We first establish a weighted estimate for $\theta_{L(t)}$ in the inner region in terms of its values on the outer annulus
${C^{-1}R(t)\leq|\mathbf x|\leq R(t)}$
(Theorem~\ref{thm angle est}).
This estimate is obtained by a blow-up argument using Liouville-type theorems for ancient solutions of the heat equation on the limiting Lawlor neck and punctured two-plane cone appearing in items~\ref{assumption, blow-up to Lawlor neck} and~\ref{assumption, blow-up to cone} of Assumption~\ref{blow-up assumptions}.
Interior gradient estimates for $\theta_{L(t)}$
(Theorem~\ref{thm deriv angle est}), together with the identity
$\mathbf H_{L(t)}=J\nabla\theta_{L(t)}$, then yield the following pointwise mean curvature bound.

\begin{theorem}\label{cor: pointwise mean curvature bound}
    Let $\{F_t:L\to L(t)\}_{t\in[0,T)}$ satisfy
    Assumption~\ref{blow-up assumptions}
    with inner scale $\lambda(t)$ and outer scale $R(t)$. Define a regularized distance function by $\rho(x,t) := \sqrt{ \lambda(t)^2 + |F_t(x)|^2 }$ for $(x,t) \in L \times [0, T)$.
    Let $\theta(x,t)$ be a globally defined Lagrangian angle,
    and let $\mathbf H(x,t)$ denote the mean curvature vector
    along $F_t$. Fix a constant $\overline\theta\in\R$.

    For every $\alpha\in(0,n-2)$ and $\beta\in\R$,
    there exist $C>1$ and $t_0\in[0,T)$ such that
    \begin{equation} \label{cor: pointwise mean curvature bound, eqn 1}
        |\mathbf H(x,t)|
        \leq C\left(
            1+
            \frac{\lambda(t)^\beta}
                 {\rho(x,t)^{1+\alpha}}
            \Theta_{\alpha,\beta,C}
        \right)
    \end{equation}
    whenever $t\in[t_0,T)$ and
    $|F_t(x)|\leq C^{-1}R(t)$, where
    \[
        \Theta_{\alpha,\beta,C}
        :=
        \sup_{t'\in[0,T)}
        \sup_{\substack{x'\in L\\
            C^{-1}R(t')\leq |F_{t'}(x')|\leq R(t')}}
        \frac{\rho(x',t')^\alpha}{\lambda(t')^\beta}
        |\theta(x',t')-\overline\theta|.
    \]
\end{theorem}

This estimate propagates Lagrangian angle bounds at the outer scale $R(t)$ to mean curvature bounds at smaller scales. 
We apply this estimate to the finite-time Lawlor neck pinches
$\{L_K(t)\subset\C^n\}_{t\in[0,T)}$ constructed in \cite{SS26II} for admissible $K\geq 3$, and to the infinite-time
Lawlor neck pinches $\{N(t)\subset\C^n/\Lambda\}_{t\in[t_0,\infty)}$ constructed in \cite{STW}. 
In each case, we verify the structural assumptions and the weighted angle estimates needed to bound the right-hand side of \eqref{cor: pointwise mean curvature bound, eqn 1}. 
In the flat-torus setting, the structural assumptions are verified using Euclidean lifts
near the singular points. Combining the resulting local estimates with bounds away from the neck yields the existence of both finite-time and infinite-time LMCF singularities with bounded mean curvature in dimensions $n \ge 3$.

\begin{corollary} \label{main cor intro}
Let $n\geq 3$.
    \begin{enumerate}
        \item (Corollary \ref{main cor finite-time sings sect})
        The mean curvature of the finite-time Lawlor neck pinches
        $\{L_K(t)\subset\mathbb{C}^{n}\}_{t\in[0,T)}$ constructed in \cite{SS26II} with curvature blow-up $\sim(T-t)^{-K/2}$ for admissible $K\geq 3$ remains
        uniformly bounded as $t\to T$.

        \item (Corollary \ref{main cor infinite-time sings sect})
        The mean curvature of the infinite-time Lawlor neck pinches
        $\{N(t)\subset\mathbb{T}^{n}\}_{t\in[t_0,\infty)}$ constructed in \cite{STW} with curvature blow-up $\sim t^{\frac{1}{n-2}}$ for admissible parameters remains
        uniformly bounded as $t\to\infty$.
    \end{enumerate}
\end{corollary}
See Corollaries \ref{main cor finite-time sings sect} and \ref{main cor infinite-time sings sect} for precise statements and notions of admissible parameters.
These examples show that, in complex dimensions $n\geq 3$, both finite-time and infinite-time Lawlor neck pinches can occur with uniformly bounded mean curvature.
Very recently, Lotay--Oliveira \cite{LO26} proved that the mean curvature becomes unbounded at nondegenerate neck pinch singularities of LMCF in compact Calabi--Yau surfaces. 

Finally, we note that the limiting two-plane cones and Lawlor necks in Assumption~\ref{blow-up assumptions} may a priori depend on the choice of subsequence. The first-named author \cite{Stolarski26} proved that, for integral Brakke flows in Euclidean space with locally uniformly bounded areas and mean curvature, a tangent flow given by a multiplicity-one regular stationary cone is unique. Theorem~\ref{cor: pointwise mean curvature bound} may therefore provide a route to proving uniqueness of such singular tangent flows for LMCF in $\mathbb C^n$, $n\geq 3$.

For finite-time singularities, a natural candidate for the outer scale is $R(t)\sim\sqrt{T-t}$. The remaining challenges are to identify an inner scale $\lambda(t)$ satisfying Assumption~\ref{blow-up assumptions} and to establish weighted Lagrangian angle bounds on the outer annulus strong enough to yield uniform mean curvature bounds. Under suitable hypotheses, Li--Sz\'ekelyhidi \cite{LS24} proved the existence of a Type II blow-up sequence converging to a Lawlor neck at singularities admitting a multiplicity-one tangent flow given by two transverse special Lagrangian planes. Their result thus provides candidate inner scales along a blow-up sequence.

\addtocontents{toc}{\protect\setcounter{tocdepth}{0}} 
\subsection*{Acknowledgment} 
The first named-author was supported by a Leverhulme Trust Early Career Fellowship (ECF-2023-182).
The second-named author is supported by the Taiwan NSTC grant (114-2115-M-008-012-MY3). 

\subsection*{Statement on AI use} AI tools were used only during the final stages of manuscript preparation to improve the exposition and ensure consistency of notation. They were not used to develop the mathematical ideas presented in this paper, nor the companion papers \cites{SS26I, SS26II}.
\addtocontents{toc}{\protect\setcounter{tocdepth}{2}}

\section{Estimates at Small Scales} \label{Section Ests at Small Scales}

This section proves the pointwise mean curvature estimate Theorem \ref{cor: pointwise mean curvature bound}, which bounds the mean curvature at small scales near a singularity by the Lagrangian angle at larger scales.
Throughout this section, we frequently invoke Assumption \ref{blow-up assumptions}, which asserts the existence of inner and outer scale functions $\lambda(t)$ and $R(t)$, respectively, and characterizations of certain blow-up limits of the Lagrangian mean curvature flow.

The comparison properties (Assumption \ref{blow-up assumptions} \eqref{scale assumptions}) of the scale functions $R, \lambda$ imply the following property that we use often in our subsequent blow-up analysis.

\begin{lem} \label{lem lambda ratio converges to 1}
    Assume $0 < T \le \infty$ and functions $R: [0,T) \to (0, r_0]$ and $\lambda: [0, T) \to (0, \infty)$ satisfy Assumption \ref{blow-up assumptions} \eqref{scale assumptions}.

    For all sequences $t_i \in [0, T)$ with $t_i \to T$ and $\rho_i > 0$ with $\frac{\rho_i}{R(t_i)} \to 0$, 
    for any $\sigma \ge 0$, 
    \begin{equation}
        \lim_{i \to \infty} \sup_{t \in [t_i - \sigma \rho_i^2, t_i]} \left| \frac{\lambda(t)}{\lambda(t_i)} - 1 \right| = 0.
    \end{equation}
\end{lem}        
\begin{proof}
    Let $\sigma \ge 0$.
    For any $t \in [t_i - \sigma \rho_i^2, t_i]$, Assumption \ref{blow-up assumptions} \eqref{scale assumptions} implies
    \begin{equation*}
        \left|\log\frac{\lambda(t)}{\lambda(t_{i})}\right|\leq\int_{t}^{t_{i}}\left|\frac{\lambda'(s)}{\lambda(s)}\right|\:ds\leq D_0R(t_{i})^{-2}(t_{i} - t)\leq D_0\sigma\frac{\rho_{i}^{2}}{R(t_{i})^{2}}\xrightarrow[i \to \infty]{} 0.
    \end{equation*}
    Exponentiating and taking a supremum in $t$ then completes the proof.
\end{proof}
     
\begin{definition}
    For a zero-Maslov Lagrangian mean curvature flow $\left\{ F_t : L \to  L(t) \right\}_{t \in [0, T)}$ with outer and inner scales $R(t), \lambda(t)$ as in Assumption \ref{blow-up assumptions}, we denote
        $$\rho(x,t) := \sqrt{ \lambda(t)^2 + |F_t(x)|^2 }   \qquad \text{for } (x,t) \in L \times [0, T).$$
\end{definition}

\begin{remark}
    The arguments in this section all proceed by a blow-up analysis around a singular point. They therefore continue to hold more generally for flows in an ambient Calabi-Yau manifold $X^n$ instead of $\C^n$.
    In this generality, the statements and proofs follow similarly except with the modulus $| \cdot |$ replaced by the ambient distance to the singular point.
\end{remark}

The following notation will be convenient for simplifying statements.
\begin{definition}[\cite{SS26I}*{Definition \ref{I-defn notation}}]
    For general nonnegative functions $f = f( x) \ge 0 $ and $g = g( x) \ge 0$, we use the notation ``$f \lesssim g$'' to mean there exists $C > 0$ such that $f(x) \le C g ( x)$ for all $ x$.
    We write ``$f \lesssim_{a,b,c} g $'' to indicate that $f \lesssim g$ and the implicit constant $C = C(a,b,c)$ depends on $a,b,c$.
    Analogous definitions apply for additional or fewer parameters $a,b,c$.
    The notation ``$f \sim g$'' means $f \lesssim g$ and $g \lesssim f$.
    For sequences $(a_i)_{i \in \mathbb N}$ and $(b_i)_{i \in \mathbb N}$ of positive real numbers, we use the notation ``$a_i \ll b_i$'' to mean $\lim_{i \to \infty} \frac{a_i}{b_i} = 0$.
\end{definition}

The following key result bounds the Lagrangian angle \emph{within} the outer scale by its values \emph{at} the outer scale.
Its proof is based on a blow-up and contradiction argument with the Liouville-type theorems from \cite{STW}*{Propositions 6.4 and 6.5}.
    
\begin{theorem} \label{thm angle est}
    Let $\left\{ F_t : L \to L(t) \right\}_{t \in [0, T)}$ be a zero-Maslov Lagrangian mean curvature flow with outer and inner scales $R(t), \lambda(t)$ satisfying Assumption \ref{blow-up assumptions}.
    Let $\theta = \theta_{L(t)}$ be the Lagrangian angle of $L(t)$ and $\overline \theta \in \R$ a constant.

    If $\alpha \in (0, n-2)$ and $\beta \in \R$, then there exists a constant $C > 1$ such that 
    \begin{equation} \label{eqn angle est}
        \sup_{t \in [0, T)} \sup_{\substack{x \in L, \\| F_t(x) | \le  R(t)}} \frac{\rho(x,t)^\alpha}{\lambda(t)^\beta } \cdot  | \theta(x,t) - \overline \theta| 
        \le C \sup_{t \in [0, T)} \sup_{\substack{x \in L, \\ C^{-1} R(t) \le | F_t(x) | \le  R(t)}} \frac{\rho(x,t)^\alpha}{\lambda(t)^\beta } \cdot  | \theta(x,t) - \overline \theta| .
    \end{equation}
\end{theorem}
\begin{proof}
    Let $\alpha \in (0, n-2)$ and $\beta \in \R$.
    Without loss of generality, we can assume the quantity on the right-hand side of \eqref{eqn angle est} is finite, that is,
        $$\sup_{t \in [0, T)} \sup_{\substack{x \in L, \\ C^{-1} R(t) \le | F_t(x) | \le  R(t)}} \frac{\rho(x,t)^\alpha}{\lambda(t)^\beta } \cdot  | \theta - \overline \theta|  < \infty \qquad \text{for all } C > 1,$$
    otherwise the conclusion is trivial.
    We can also assume this quantity is non-zero for all $C > 1$ by unique continuation for $\partial_t (\theta- \overline \theta) = \Delta_{L(t)} (\theta - \overline \theta)$.

    For the sake of contradiction, suppose \eqref{eqn angle est} were false.
    Then, for any $C_i = i \in \mathbb N$, 
        $$\sup_{t \in [0, T)} \sup_{\substack{x \in L, \\| F_t(x) | \le  R(t)}} \frac{\rho(x,t)^\alpha}{\lambda(t)^\beta } \cdot  | \theta - \overline \theta| 
        > i \sup_{t \in [0, T)} \sup_{\substack{x \in L, \\ \frac1i R(t) \le | F_t(x) | \le  R(t)}} \frac{\rho(x,t)^\alpha}{\lambda(t)^\beta } \cdot  | \theta - \overline \theta|.$$
    It follows that there is an increasing sequence of times $T_i \in [0, T)$ such that 
    \begin{equation*}
        M_i := \sup_{0 \le t \le T_i} \sup_{\substack{x \in L, \\| F_t(x) | \le  R(t)}} \frac{\rho(x,t)^\alpha}{\lambda(t)^\beta } \cdot  | \theta - \overline \theta| 
        > i \sup_{t \in [0, T)} \sup_{\substack{x \in L, \\ \frac1i R(t) \le | F_t(x) | \le  R(t)}} \frac{\rho(x,t)^\alpha}{\lambda(t)^\beta } \cdot  | \theta - \overline \theta| \xrightarrow[i \to \infty]{} \infty.
    \end{equation*}
    By continuity of $\frac{\rho(x,t)^\alpha}{\lambda(t)^\beta } \cdot  | \theta(x,t) - \overline \theta| $ and compactness of the set where $0 \le t \le T_i$, $x \in L$, and $|F_t(x)| \le R(t)$, 
    it follows that there exists a sequence of spacetime points $(x_i, t_i)$ such that $0 \le t_i \le T_i$, $x_i \in L$, $|F_{t_i} (x_i)| \le R(t_i)$ and
    \begin{equation*}
        M_i = \frac{\rho(x_i,t_i)^\alpha}{\lambda(t_i)^\beta } \cdot  | \theta (x_i, t_i) - \overline \theta| .
    \end{equation*}
    In summary, combining the above results implies
    \begin{multline} \label{proof angle est, eqn 1}
        M_i
        = \frac{\rho(x_i,t_i)^\alpha}{\lambda(t_i)^\beta } \cdot  | \theta (x_i, t_i) - \overline \theta| 
        =\sup_{0 \le t \le t_i} \sup_{\substack{x \in L, \\| F_t(x) | \le  R(t)}} \frac{\rho(x,t)^\alpha}{\lambda(t)^\beta } \cdot  | \theta - \overline \theta| \\
        > i \sup_{t \in [0, T)} \sup_{\substack{x \in L, \\ \frac1i R(t) \le | F_t(x) | \le  R(t)}} \frac{\rho(x,t)^\alpha}{\lambda(t)^\beta } \cdot  | \theta - \overline \theta| \xrightarrow[i \to \infty]{} \infty.
    \end{multline}
    It must then be the case that $t_i \to T$, otherwise $M_i$ would be bounded. 
    Additionally, \eqref{proof angle est, eqn 1} implies $|F_{t_i}(x_i) | \le \frac1i R(t_i)$ and so $\lim_{i \to \infty} \frac{|F_{t_i}(x_i)|}{R(t_i)} = 0$.

    There are two cases to consider based on the behavior of the spacetime sequence $(x_i, t_i)$: either $|F_{t_i}(x_i)| \lesssim \lambda(t_i)$ along a subsequence or $|F_{t_i}(x_i) | \gg \lambda(t_i)$ along a subsequence.
    
    \underline{\textit{Case 1: $| F_{t_i}(x_i) | \lesssim \lambda(t_i)$ along a subsequence.}}
    
    Assume first that there exists a subsequence (still denoted by index $i$) and a constant $C_0' > 0$ such that $|F_{t_i}(x_i) | \le C_0' \lambda(t_i)$ for all $i$.
    To simplify the notation, denote $\lambda_i := \lambda(t_i)$.
    Consider the sequence of parabolically rescaled Lagrangian mean curvature flows
        $$\tilde L_i (\tau) := \frac1{\lambda_i} L ( t_i + \tau \lambda_i^2), \qquad \tau \in \left[ - \frac{t_i}{\lambda_i^2}, 0 \right]$$
    given by the image of the immersions
        $$\tilde F_{i, \tau} : L \to \tilde L_i(\tau) \subset  \C^n, \qquad \tilde F_{i, \tau}(x) = \frac1{\lambda_i} F_{t_i + \tau \lambda_i^2} (x).$$
    By Assumption \ref{blow-up assumptions} \eqref{assumption, blow-up to Lawlor neck}, $\tilde L_i (\tau)$ converges in $C^\infty_{\loc}( \C^n \times (-\infty, 0]) $ to a Lawlor neck $\overline L$ as $i \to \infty$.

    Denote the Lagrangian angles of the rescaled flows by $\tilde \theta_i := \theta_{\tilde L_i}$.    
    Recall $\partial_\tau ( \tilde \theta_i - \overline \theta) = \Delta_{\tilde L_i(\tau)} ( \tilde \theta_i - \overline \theta)$.
    Define functions 
        $$\tilde u_i : L \times \left[ - \frac{t_i}{\lambda_i^2} , 0\right]  
        \to \R \qquad \text{by} \qquad 
        \tilde u_i(x, \tau) := \frac{\lambda_i^\alpha}{M_i \lambda_i^{ \beta}} \cdot ( \tilde \theta_i (x, \tau) - \overline \theta)  .$$ 
    Then
    \begin{equation} \label{proof angle est, eqn 1.5}
        \partial_\tau \tilde u_i = \Delta_{\tilde L_i(\tau)} \tilde u_i \qquad \text{on } L \times \left[-\frac{t_i}{\lambda_i^2} , 0\right].
    \end{equation}
    Using \eqref{proof angle est, eqn 1}, it follows that, for all $i \in \N$ and $(x,\tau) \in L \times \left[-\frac{t_i}{\lambda_i^2}, 0 \right]$,
    \begin{gather} \label{proof angle est, eqn 2} \begin{aligned}
        &|\tilde u_i ( x , \tau)| \\
        ={}& \frac{\rho(x,t)^{-\alpha}}{\lambda(t)^{-\beta}} 
        \frac{\rho(x,t)^\alpha}{\lambda(t)^\beta}  \frac{\lambda_i^\alpha}{M_i \lambda_i^{ \beta}}| \theta( x, t) - \overline \theta| 
        && ( \text{if }  t = t_i + \lambda_i^2\tau  )\\
        \le{}& \frac{\rho(x, t)^{-\alpha} \lambda_i^\alpha}{\lambda(t)^{-\beta} \lambda_i^{ \beta}} 
        && (\text{if } |F_t(x)| \le R(t) \text{ and } t \le t_i )\\
        ={}& \left(\frac{\lambda^2(t_i + \tau \lambda_i^2)}{\lambda_i^2} + |\tilde F_{i, \tau}(x)|^2 \right)^{-\alpha/2} \cdot \left( \frac{\lambda(t_i + \tau \lambda_i^2)}{\lambda_i} \right)^\beta
    \end{aligned} \end{gather}       
    where $\rho(x,t) = \sqrt{ \lambda(t)^2 + |F_t(x)|^2} $ above is still defined with respect to the flow $L(t)$.
    Since $R(t)$ is non-increasing, we deduce from \eqref{proof angle est, eqn 2} that, for all $i \in \N$, 
    \begin{multline} \label{proof angle est, eqn 3}
        |\tilde u_i (x,\tau) | \le \left(\frac{\lambda^2(t_i + \tau \lambda_i^2)}{\lambda_i^2} + |\tilde F_{i, \tau}(x)|^2 \right)^{-\alpha/2} \cdot \left( \frac{\lambda(t_i + \tau \lambda_i^2)}{\lambda(t_i)} \right)^\beta
        \\
        \forall (x, \tau) \in  L \times \left[ - \frac{t_i}{\lambda_i^2}, 0 \right] \text{ with } |\tilde F_{i, \tau}(x)| \le \frac{R(t_i)}{\lambda_i}.
    \end{multline}
    
    By Assumption \ref{blow-up assumptions} \eqref{scale assumption 1}, $\frac{\lambda_i}{R(t_i)}= \frac{\lambda(t_i)}{R(t_i)} \to 0$ as $i \to \infty$.
    Lemma \ref{lem lambda ratio converges to 1} with $\rho_i = \lambda_i$ then implies
    \begin{equation} \label{proof angle est, eqn 4}
        \frac{\lambda_i( t_i + \tau \lambda_i^2)}{\lambda_i^2} \xrightarrow[i \to \infty]{C^0_{\loc}( (-\infty, 0])} 1.
    \end{equation}

    Since $\tilde L_i(\tau) \xrightarrow[i \to \infty]{C^\infty_{\loc}(\C^n \times (-\infty, 0])} \overline L$,
    for any point $(\mathbf x, \tau) \in \overline L \times (-\infty, 0] \subset \C^n \times (-\infty, 0]$, 
    there exists a bounded spacetime neighborhood $\Omega  \subset \C^n \times (-\infty, 0]$ of $(\mathbf x, \tau)$ and $I \in \N$
    such that, for $i > I$, the flows $\tilde L_i(\tau)$ in $\Omega$ can be written as normal graphs over $(\overline L \times \R ) \cap \Omega$ and these functions converge to $0$ in $C^\infty((\overline L \times \R ) \cap \Omega)$ as $i \to \infty$.
    Through this graphical representation, the functions $\tilde u_i$ can be regarded as functions on $(\overline L \times \R) \cap \Omega$ for $i > I$.
    By increasing $I$ and invoking \eqref{proof angle est, eqn 3} and \eqref{proof angle est, eqn 4}, these functions $\tilde u_i$ have $C^0((\overline L \times \R) \cap \Omega)$-bounds that are uniform in $i >I$.
    The $\tilde u_i$ solve the heat equation \eqref{proof angle est, eqn 1.5}, which, by the convergence $\tilde L_i(\tau) \to \overline L$, is a linear parabolic PDE whose coefficients satisfy $C^\infty((\overline L \times \R) \cap \Omega)$-bounds that are uniform in $i > I$.
    Standard interior estimates for linear parabolic PDEs therefore bootstrap the $C^0((\overline L \times \R) \cap \Omega)$-bounds on $\tilde u_i$ to $C^\infty((\overline L \times \R) \cap \Omega')$-bounds on $\tilde u_i$ in some smaller spacetime neighborhood $\Omega'$ of $(\mathbf x, \tau)$, which again are uniform in $i > I$.
    
    Note that, by letting $(\mathbf x, \tau) \in \overline L \times (-\infty, 0]$ vary, the associated neighborhoods $\Omega'$ cover $\overline L \times (-\infty, 0]$.
    By extracting a countable subcover and using the Arzel{\`a}-Ascoli Theorem with a diagonalization argument, it follows that we can pass to a further subsequence of the $\tilde u_i$ that converges in $C^\infty_{\loc}(\overline L \times (-\infty,0])$ to some smooth function $\tilde u_\infty : \overline L \times (-\infty, 0] \to \R$.
    Taking the limit $i\to \infty$ in equations \eqref{proof angle est, eqn 1.5}, \eqref{proof angle est, eqn 3}, and \eqref{proof angle est, eqn 4} then implies 
    \begin{equation} \label{proof angle est, eqn 5}
        \partial_\tau \tilde u_\infty = \Delta_{\overline L} \tilde u_\infty 
        \quad \text{and} \quad |\tilde u_\infty(\mathbf x, \tau)| \le \left( 1 + |\mathbf x|^2 \right)^{-\alpha/2} \qquad \text{for all } (\mathbf x, \tau) \in \overline L \times (-\infty, 0].
    \end{equation}

    Because $\alpha \in (0, n-2)$, the Liouville-type Theorem \cite{STW}*{Proposition~6.4} gives
    \begin{equation} \label{proof angle est, eqn 6}
        \tilde u_\infty \equiv 0 \qquad \text{on } \overline L \times (-\infty, 0].
    \end{equation}

    On the other hand, \eqref{proof angle est, eqn 1} implies the inequality in \eqref{proof angle est, eqn 2} becomes an equality at $( x_i, \tau = 0)$, and so
    \begin{equation} \label{proof angle est, eqn 7}
        |\tilde u_i (x_i, 0) | = \left(1 + \frac{|F_{t_i}(x_i)|^2}{\lambda_i^2}
         \right)^{-\alpha/2}  \qquad \forall i \in \N.
    \end{equation}
    Since $| F_{t_i}(x_i) | \lesssim \lambda_i$, $\frac{F_{t_i}(x_i)}{\lambda_i}$ subsequentially converges to some $\mathbf x_\infty \in \overline L$. 
    Equation \eqref{proof angle est, eqn 7} and the $C^\infty_{\loc}(\overline L \times (-\infty, 0])$-convergence of $\tilde u_i$ to $\tilde u_\infty$ therefore implies
    \begin{equation} \label{proof angle est, eqn 8}
        |\tilde u_\infty(\mathbf x_\infty, 0)| = \left( 1 + |\mathbf x_\infty|^2 \right)^{-\alpha /2} \ne 0.
    \end{equation}
    This however contradicts \eqref{proof angle est, eqn 6}.

    \underline{\textit{Case 2: $\lambda(t_i) \ll | F_{t_i}(x_i) | \ll R(t_i)$ along a subsequence.}}
    
    Assume now that there exists a subsequence (still denoted by index $i$) such that $|F_{t_i}(x_i)| \gg \lambda_i$, that is, $\liminf_{i \to \infty} \frac{|F_{t_i}(x_i)|}{\lambda_i} = \infty$.
    Recall additionally that $\lim_{i \to \infty} \frac{|F_{t_i}(x_i)|}{R(t_i)} = 0$.
    The overall argument in this case is similar to Case 1, except we rescale the flows by $|F_{t_i}(x_i)|$ instead of $\lambda_i = \lambda(t_i)$.
    We detail the argument nonetheless.
    Consider the sequence of parabolically rescaled Lagrangian mean curvature flows
        $$\tilde L_i (\tau) := \frac1{|F_{t_i}(x_i)|} L ( t_i + \tau |F_{t_i}(x_i)|^2), \qquad \tau \in \left[ - \frac{t_i}{|F_{t_i}(x_i)|^2}, 0 \right]$$
    given by the image of the immersions   
        $$\tilde F_{i, \tau} : L \to \tilde L_i(\tau) \subset \C^n,\qquad \tilde F_{i, \tau}(x) = \frac1{|F_{t_i}(x_i)|} F_{t_i + \tau |F_{t_i}(x_i)|^2} ( x)$$
    By Assumption \ref{blow-up assumptions} \eqref{assumption, blow-up to cone}, $\tilde L_i (\tau)$ converges in $C^\infty_{\loc}( \C^n \setminus \{ \mathbf 0\}  \times (-\infty, 0])$ to a special Lagrangian cone $\Cone = \mathcal P_0 \cup \mathcal P_1$ consisting of two transversely intersecting planes as $i \to \infty$.

    Let $\tilde \theta_i := \theta_{\tilde L_i}$ denote the Lagrangian angle of the rescaled flow $\tilde L_i(\tau)$, and 
    define functions
    \begin{equation*}
        \tilde u_i : L \times \left[ - \frac{t_i}{|F_{t_i}(x_i)|^2} , 0\right] \to \R \quad \text{given by} \quad  \tilde u_i := \frac{|F_{t_i}(x_i)|^\alpha}{M_i \lambda_i^\beta} \cdot \left( \tilde \theta_i - \overline \theta \right).
    \end{equation*}
    Then
    \begin{equation} \label{proof angle est, eqn 9}
        \partial_\tau \tilde u_i = \Delta_{\tilde L_i(\tau)} \tilde u_i \qquad \text{on } L \times \left[ -\frac{t_i}{|F_{t_i}(x_i)|^2}, 0 \right],
    \end{equation}
    and \eqref{proof angle est, eqn 1} implies, for all $i\in \N$ and $(x, \tau) \in L \times \left[ - \frac{t_i}{|F_{t_i}(x_i)|^2}, 0 \right]$,
    \begin{gather} \label{proof angle est, eqn 10} \begin{aligned}
        |\tilde u_i (x, \tau) |
        ={}& \frac{\rho(x,t)^{-\alpha}}{\lambda(t)^{-\beta}} 
        \frac{\rho(x,t)^\alpha}{\lambda(t)^\beta}  \frac{|F_{t_i}(x_i)|^\alpha}{M_i \lambda_i^{ \beta}}| \theta( x, t) - \overline \theta| 
        && ( \text{if } t = t_i + |F_{t_i}(x_i)|^2\tau  )\\
        \le{}& \frac{\rho(x,t)^{-\alpha} |F_{t_i}(x_i)|^\alpha}{\lambda(t)^{-\beta} \lambda_i^{ \beta}} 
        && (\text{if } |F_{t}(x)| \le R(t) \text{ and } t \le t_i )\\
        ={}& \left(\frac{\lambda^2(t_i + \tau |F_{t_i}(x_i)|^2)}{|F_{t_i}(x_i)|^2} + |\tilde F_{i, \tau}(x)|^2 \right)^{-\alpha/2} \\& \qquad \cdot \left( \frac{\lambda(t_i + \tau |F_{t_i}(x_i)|^2)}{\lambda_i} \right)^\beta.
    \end{aligned} \end{gather}
    Thus, for all $i \in \N$, 
    \begin{multline} \label{proof angle est, eqn 11}
        |\tilde u_i (x, \tau) | \le \left(\frac{\lambda^2(t_i + \tau |F_{t_i}(x_i)|^2)}{\lambda_i^2} \cdot \frac{\lambda_i^2}{|F_{t_i}(x_i)|^2} + |\tilde F_{i, \tau}(x)|^2 \right)^{-\alpha/2} \cdot \left( \frac{\lambda(t_i + \tau |F_{t_i}(x_i)|^2)}{\lambda_i} \right)^\beta \\
        \forall (x, \tau) \in L \times \left[ - \frac{t_i}{|F_{t_i}(x_i)|^2} , 0\right] \text{ with } |\tilde F_{i, \tau}(x)| \le \frac{R(t_i)}{|F_{t_i}(x_i)|}.
    \end{multline} 
    Recall $\frac{\lambda_i}{|F_{t_i}(x_i)|} \to 0$ and $\frac{|F_{t_i}(x_i)|}{R(t_i)} \to 0$ as $i \to \infty$.
    Lemma \ref{lem lambda ratio converges to 1} with $\rho_i = |F_{t_i}(x_i)|$ then implies
    \begin{equation} \label{proof angle est, eqn 12}
        \frac{\lambda_i( t_i + \tau |F_{t_i}(x_i)|^2)}{\lambda_i^2} \xrightarrow[i \to \infty]{C^0_{\loc}( (-\infty, 0])} 1.
    \end{equation}
    Recall $\tilde L_i (\tau) \xrightarrow[i \to \infty]{C^\infty_{\loc}(\C^n \setminus \{ \mathbf 0\} \times (-\infty , 0])} \Cone$.
    By similar logic as in Case 1 above, the functions $\tilde u_i$ can be regarded as a sequence of functions defined on a sequence of spacetime domains exhausting $\Cone \setminus \{ \mathbf 0 \} \times (-\infty, 0]$.
    Combined with \eqref{proof angle est, eqn 11}, it follows that the functions $\tilde u_i$ satisfy $C^0_{\loc}(\Cone \setminus \{ \mathbf 0 \}  \times (-\infty, 0])$-estimates that are uniform in $i$.
    Since the $\tilde u_i$ solve the heat equation \eqref{proof angle est, eqn 9} and $\tilde L_i (\tau) \xrightarrow[i \to \infty]{C^\infty_{\loc}(\C^n \setminus \{ \mathbf 0\} \times (-\infty , 0])} \Cone$, interior estimates for linear parabolic PDEs bootstrap these $C^0_{\loc}$-estimates to $C^\infty_{\loc}$-estimates for $\tilde u_i$ which are uniform in $i$.
    By the Arzel{\`a}-Ascoli Theorem, we can therefore pass to a subsequence of the $\tilde u_i$ that $C^\infty_{\loc}(\Cone \setminus \{ \mathbf 0 \} \times (-\infty, 0])$-converges to some smooth function $\tilde u_\infty : \Cone \setminus \{ \mathbf 0 \} \times (-\infty, 0] \to \R$.
    Taking the limit $i\to \infty$ in equations \eqref{proof angle est, eqn 9}, \eqref{proof angle est, eqn 11}, and \eqref{proof angle est, eqn 12} then implies 
    \begin{equation} \label{proof angle est, eqn 13}
        \partial_\tau \tilde u_\infty = \Delta_{\Cone} \tilde u_\infty 
        \quad \text{and} \quad |\tilde u_\infty(\mathbf x, \tau)| \le |\mathbf x|^{-\alpha} \qquad \text{on } \Cone \setminus \{ \mathbf 0 \}  \times (-\infty, 0].
    \end{equation}
    Because $\alpha \in (0, n-2)$, the Liouville-type Theorem \cite{STW}*{Proposition~6.5} gives
    \begin{equation} \label{proof angle est, eqn 14}
        \tilde u_\infty \equiv 0 \qquad \text{on } \Cone \setminus \{ \mathbf 0 \} \times (-\infty, 0].
    \end{equation}

    On the other hand, \eqref{proof angle est, eqn 1} implies the inequality in \eqref{proof angle est, eqn 10} becomes an equality at $(x_i, \tau = 0)$, and so
    \begin{equation} \label{proof angle est, eqn 15}
        |\tilde u_i (x_i, 0) | = \left(\frac{\lambda_i^2}{|F_{t_i}(x_i)|^2} + 1
         \right)^{-\alpha/2}  \qquad \forall i \in \N.
    \end{equation}
    The sequence $\tilde F_{i, 0}(x_i) = \frac{F_{t_i}(x_i)}{|F_{t_i}(x_i)|}$ subsequentially converges to some $\mathbf x_\infty \in \Cone$ with $|\mathbf x_\infty| = 1$. 
    This convergence, equation \eqref{proof angle est, eqn 15}, and the $C^\infty_{\loc}$-convergence of $\tilde u_i$ to $\tilde u_\infty$ therefore implies
    \begin{equation} \label{proof angle est, eqn 16}
        |\tilde u_\infty(\mathbf x_\infty, 0)| = \left( 0 + 1 \right)^{-\alpha /2} \ne 0.
    \end{equation}
    This however contradicts \eqref{proof angle est, eqn 14}.
        
    Because both cases lead to a contradiction, the proof of the theorem is complete.
\end{proof}

For the next result, it will be convenient to introduce the following notation.
\begin{definition}
    For $(\mathbf x_0, t_0) \in \C^n \times \R$ and $r > 0$, denote the closed backwards parabolic neighborhood
        $$P(\mathbf x_0, t_0; r) := \overline B(\mathbf x_0, r) \times [t_0 - r^2, t_0] \subseteq \C^n \times \R.$$
    If $\{ F_t : L \to L(t) \}_{t \in [0, T)}$ is a Lagrangian mean curvature flow in $\C^n$, the restriction of $P (\mathbf x_0, t_0; r)$ to the spacetime track of the flow $L(t)$ is denoted
        $$P_L(\mathbf x_0, t_0; r) := \{ (\mathbf x, t) \in \C^n \times [0, T) \, \colon \,  \, \mathbf x \in L(t), \, | \mathbf x - \mathbf x_0 | \le r, \text{ and } t \in [t_0 - r^2, t_0]  \}.$$
    For $x_0 \in L$, we let $P_L (x_0, t_0; r)$ denote the pre-image of $P_L ( F_{t_0}(x_0), t_0; r) $ under $F_t$, that is
        $$P_L (x_0, t_0; r) := \left\{ (x, t) \in L \times [0, T) \colon |F_t(x) - F_{t_0}(x_0)| \le r \text{ and } t \in [t_0 - r^2, t_0]\right \}.$$
\end{definition}

By leveraging interior regularity results for parabolic PDEs with a blow-up argument, the next result transfers estimates on the Lagrangian angle to estimates on its derivative.
\begin{theorem} \label{thm deriv angle est}
    Let $\left\{ F_t : L \to L(t) \right\}_{t \in [0, T)}$ be a zero-Maslov Lagrangian mean curvature flow with outer and inner scales $R(t)$ and $\lambda(t)$, respectively, satisfying Assumption \ref{blow-up assumptions}.
    Let $\theta = \theta_{L(t)}$ be the Lagrangian angle of $L(t)$ and $\overline \theta \in \R$ a constant.

    Then there exists a constant $C>1$ such that
    \begin{equation}
        |d \theta|_{L(t)} (x, t) \le C \left( 1 + \frac1{\rho(x,t)} \cdot \sup_{P_L\left(x, t; \frac12 \rho(x,t) \right)} \left|\theta - \overline \theta \right| \right) 
    \end{equation}
    for all $(x, t)  \in L \times [0, T)$ with  $| F_t(x) | \le \frac1C R(t)$.
\end{theorem}

\begin{proof}
    Suppose otherwise for the sake of contradiction.
    Then for each $i \in \mathbb N$ there exists $(x_i, t_i) \in L \times [0, T)$ such that $|F_{t_i}(x_i)| \le \frac1i R(t_i)$,
    and, denoting $\rho_i := \rho(x_i, t_i) = \sqrt{ \lambda(t_i) + |F_{t_i}(x_i)|^2}$, 
    \begin{equation} \label{proof deriv angle est, eqn 1}
        | d \theta|_{L(t_i)} (x_i, t_i) 
        > i \left( 1 + \frac1{\rho_i} \sup_{P_L \left(x_i, t_i; \frac12 \rho_i \right)} |\theta - \overline \theta | \right) \ge i.
    \end{equation}
    It must then be the case that $t_i \to T$ as $i \to \infty$.
    Thus, $|F_{t_i}(x_i) | \le \frac 1i R(t_i) \to 0$, $\lambda_i := \lambda(t_i) \to 0$, and $\rho_i \to 0$ as $i \to \infty$.
    In particular, $|F_{t_i}(x_i)| \ll R(t_i)$.

    Consider the parabolically rescaled Lagrangian mean curvature flows
        $$\tilde L_i(\tau) := \frac1{\rho_i} L( t_i + \tau \rho_i^2) , \qquad \tau \in \left[ - \frac{t_i}{\rho_i^2}, 0 \right]$$
    parametrized by
        $$\tilde F_{i, \tau} : L \to \tilde L_i(\tau) , \qquad \tilde F_{i, \tau}(x) = \frac1{\rho_i} F_{t_i + \tau \rho_i^2}(x).$$
    Let $\tilde \theta_i := \theta_{\tilde L_i}$ denote the Lagrangian angle of these rescaled flows.
    Note
    \begin{equation} \label{proof deriv angle est, eqn 2}
        \partial_\tau \left( \tilde \theta_i - \overline \theta \right) = \Delta_{\tilde L_i(\tau)} \left( \tilde \theta_i - \overline \theta \right)  \qquad \text{on } L \times \left[ - \frac{t_i}{\rho_i^2}, 0 \right],
    \end{equation}
    Under this rescaling, equation \eqref{proof deriv angle est, eqn 1} becomes, for all $i \in \N$,
    \begin{gather} \label{proof deriv angle est, eqn 4} \begin{aligned}
        |d \tilde \theta_i|_{\tilde L_i(0)} (x_i, 0) 
        ={}& \rho_i\cdot | d \theta |_{L(t_i)} (x_i, t_i) \\
        >{}& \rho_i \cdot i \left( 1 + \frac1{\rho_i} \sup_{P_L \left(x_i, t_i; \frac12 \rho_i \right)} |\theta - \overline \theta | \right) \\
        ={}& \rho_i \cdot i \left( 1 + \frac1{\rho_i} \sup_{ P_{\tilde L_i} \left(x_i, 0; \frac12 \right)} |\tilde \theta_i - \overline \theta | \right) \\
        >{}& i\sup_{ P_{\tilde L_i} \left( x_i, 0; \frac12 \right)} | \tilde \theta - \overline \theta | .
    \end{aligned} \end{gather}
    
    Recall $\frac{|F_{t_i}(x_i) |}{R(t_i)}\le \frac1i \to 0$ as $i \to\infty$. There are thus two cases to consider based on the behavior the sequence $(x_i, t_i)$: either $|F_{t_i}(x_i)| \lesssim \lambda_i = \lambda(t_i)$ along a subsequence or $\lambda_i \ll |F_{t_i}(x_i) | \ll R(t_i)$ along a subsequence.

    \underline{\textit{Case 1: $| F_{t_i}(x_i) | \lesssim \lambda(t_i)$ along a subsequence.}}
    
    Assume first that there exists a subsequence (still denoted by index $i$) and a constant $C_0' > 0$ such that $|F_{t_i}(x_i) | \le C_0' \lambda_i$ for all $i$.
    Then, for all $i \in \N$, 
        $$\lambda_i \le \rho_i = \sqrt{ \lambda_i^2 + |F_{t_i}(x_i)|^2 } \le \lambda_i \sqrt{1 + C_0'^2} .$$
    By passing to a further subsequence if necessary, we can ensure the limit $\lim_{i \to \infty} \frac{\rho_i}{\lambda_i}$ exists and lies in $\left[ 1 , \sqrt{1+ C_0'^2} \right]$.
    Again passing to a further subsequence, Assumption \ref{blow-up assumptions} \eqref{assumption, blow-up to Lawlor neck} then implies the rescaled flows $\tilde L_i(\tau)$ converge in $C^\infty_{\loc}(\C^n \times (-\infty, 0])$ to a Lawlor neck $\overline L$ as $i \to \infty$.
    Additionally, $|\tilde F_{i, 0}(x_i)| = \frac{|F_{t_i}(x_i)|}{\rho_i} \le C_0'$, and so, after passing to a subsequence, the $\tilde F_{i, 0}(x_i)$ converge to some $\mathbf x_\infty \in \overline L$ with $|\mathbf x_\infty| \le C_0'$.

    By \eqref{proof deriv angle est, eqn 2}, the functions $\tilde \theta_i - \overline \theta$ satisfy the heat equation
    \begin{equation} \label{proof deriv angle est, eqn 5} 
        \partial_\tau \left( \tilde \theta_i - \overline \theta \right) = \Delta_{\tilde L_i(\tau)} \left( \tilde \theta_i - \overline \theta \right)  \qquad \text{on } P_{\tilde L_i} \left(x_i , 0; 1 \right) \subseteq L \times \left[ - \frac{t_i}{\rho_i^2}, 0 \right],
    \end{equation}
    Because the flows $\tilde L_i(\tau)$ smoothly converge to $\overline L$ in $P(\mathbf x_\infty, 0; 2)$, the linear parabolic PDE \eqref{proof deriv angle est, eqn 5} in coordinates has coefficients which are bounded in $C^3(P_{\tilde L_i}( x_i, 0; 1))$ uniformly in $i$.
    Standard interior derivative estimates for linear parabolic PDEs therefore imply there exists some $C_1' > 1$ \emph{independent of $i \in \N$} such that 
    \begin{equation} \label{proof deriv angle est, eqn 6} 
        |d \tilde \theta_i |_{\tilde L_i(0)} (x_i, 0) \le C_1' \sup_{P_{\tilde L_i}\left( x_i, 0; \frac12 \right)} \left| \tilde \theta_i - \overline \theta \right| \qquad \text{for all } i \in \N.
    \end{equation}
    This estimate \eqref{proof deriv angle est, eqn 6}, however, contradicts the inequality \eqref{proof deriv angle est, eqn 4} for $i$ sufficiently large.

    \underline{\textit{Case 2: $\lambda(t_i) \ll | F_{t_i}(x_i) | \ll R(t_i)$ along a subsequence.}}
    
    Assume in this case that there exists a subsequence (still denoted by index $i$) such that $\liminf_{i \to \infty} \frac{|F_{t_i}(x_i)|}{\lambda_i} = \infty$.
    The argument here is similar to the previous case, but we provide the details nonetheless.
    In this case,
    $$1 \le \liminf_{i \to \infty} \frac{\rho_i}{|F_{t_i}(x_i)|} \le \limsup_{i \to \infty} \frac{\rho_i}{|F_{t_i}(x_i)|}= \limsup_{i \to \infty} \sqrt{ \frac{\lambda_i^2}{|F_{t_i}(x_i)|^2} + 1 }=1.$$
    Thus, $\lim_{i \to \infty} \frac{\rho_i}{|F_{t_i}(x_i)|} = 1$ and $\lambda_i \ll \rho_i \sim |F_{t_i}(x_i)| \ll R(t_i)$.
    Assumption \ref{blow-up assumptions} \eqref{assumption, blow-up to cone} therefore implies that, after passing to a further subsequence, the rescaled flows $\tilde L_i(\tau)$ converge in $C^\infty_{\loc}(\C^n \setminus \{ \mathbf 0 \} \times (-\infty, 0])$ to a special Lagrangian cone $\Cone = \mathcal P_0 \cup \mathcal P_1$ consisting of two transversely intersecting planes.
    We can additionally deduce the $\tilde F_{i, 0}(x_i) = \frac{F_{t_i}(x_i)}{\rho_i}$ converge to some $\mathbf x_\infty \in \Cone$ with $|\mathbf x_\infty|=1$ after perhaps passing to a further subsequence.

    By \eqref{proof deriv angle est, eqn 2}, the functions $\tilde \theta_i - \overline \theta$ satisfy the heat equation
    \begin{equation} \label{proof deriv angle est, eqn 8} 
        \partial_\tau \left( \tilde \theta_i - \overline \theta \right) = \Delta_{\tilde L_i(\tau)} \left( \tilde \theta_i - \overline \theta \right)  \qquad \text{on } P_{\tilde L_i} \left(x_i , 0; \frac34 \right) \subseteq L \times \left[ - \frac{t_i}{\rho_i^2}, 0 \right],
    \end{equation}
    Note $P\left(\mathbf x_\infty, 0; \frac78 \right) \subset \C^n \setminus \{ \mathbf 0\} \times (-\infty, 0]$ since $|\mathbf x_\infty| =1$.
    Because the flows $\tilde L_i(\tau)$ smoothly converge in to $\Cone $ in $P\left(\mathbf x_\infty, 0; \frac78\right)$, the linear parabolic PDE \eqref{proof deriv angle est, eqn 8} in coordinates has coefficients which are bounded in $C^3\left(P_{\tilde L_i}\left( x_i, 0; \frac34\right) \right)$ uniformly in $i$.
    Standard interior derivative estimates for linear parabolic PDEs therefore imply there exists some $C_2' > 1$ \emph{independent of $i \in \N$} such that 
    \begin{equation} \label{proof deriv angle est, eqn 9} 
        |d \tilde \theta_i |_{\tilde L_i(0)} (x_i, 0) \le C_2' \sup_{P_{\tilde L_i}\left( x_i, 0; \frac12 \right)} \left| \tilde \theta_i - \overline \theta \right| \qquad \text{for all } i \in \N.
    \end{equation}
    This estimate \eqref{proof deriv angle est, eqn 9}, however, contradicts the inequality \eqref{proof deriv angle est, eqn 4} for $i$ sufficiently large.

    Because both cases lead to a contradiction, the proof of the theorem is complete.
\end{proof}

The proof of the mean curvature bound Theorem \ref{cor: pointwise mean curvature bound} also requires the following comparability of weights.
\begin{lemma}\label{lem: compare weights}
    Let $\left\{ F_t : L \to L(t) \right\}_{t \in [0, T)}$ be a zero-Maslov Lagrangian mean curvature flow with outer and inner scales $R(t)$ and $\lambda(t)$, respectively, satisfying Assumption \ref{blow-up assumptions}.

    For any constant $C > 2$ and any constant $D_0>0$ satisfying \eqref{eqn scale assumption 2},
    there exists $t_{0} = t_{0}(C)\in[0, T)$ and $K = K(D_{0})>0$ such that if $(x,t) \in L \times [t_{0}, T)$ with $|F_t(x)|\leq C^{-1}R(t)$, then for every $(x', t')\in P_{L}\left(x, t, \frac{1}{2}\rho(x,t) \right)$,
    \begin{align}
        &K^{-1}\lambda(t)\leq\lambda(t')\leq K\lambda(t),\label{eq: compare lambda}\\
        &K^{-1}\rho(x,t)\leq\rho(x',t')\leq K\rho(x,t),\;\;\mbox{and}\label{eq: compare rho}\\
        &|F_{t'}(x')|\leq R(t').\label{eq: compare global radius}
    \end{align}
\end{lemma}

\begin{proof}
    By Assumption \ref{blow-up assumptions} (\ref{scale assumption 1}), there is $t_{0}\in[0, T)$ such that $\lambda(t)\leq \frac{R(t)}{2}$ if $t\geq t_{0}$. Then for $ t \in [t_0, T)$ and $t'\in[t-\frac{1}{4}\rho^2(x,t), t]$,
    \begin{align*}
        t - t'\leq\frac{1}{4}\rho^{2}(x,t) = \frac{1}{4} \left(\lambda(t)^{2} + |F_t(x)|^{2} \right)\leq \left(\frac{1}{16} + \frac{1}{4C^{2}}\right)R(t)^{2}\leq\frac{1}{2}R(t)^{2}.
    \end{align*}
    Integrating \eqref{eqn scale assumption 2} and using the monotonicity of $R$ yields
    \begin{align*}
        \left|\log\frac{\lambda(t)}{\lambda(t')}\right|\leq\int_{t'}^{t}\frac{D_{0}}{R(s)^{2}}\:ds\leq\frac{D_{0}}{R(t)^{2}}\frac{1}{2}R(t)^{2}\ = \frac{D_{0}}{2}.
    \end{align*}
    This implies (\ref{eq: compare lambda}).

    We now prove (\ref{eq: compare rho}). For $(x', t')\in P_{L}\left(x, t, \frac{1}{2}\rho_{\lambda(t)}(x) \right)$ we have 
    \begin{align*}
        |F_t(x)| - \frac{1}{2}\rho(x,t)\leq|F_{t'}(x')|\leq |F_t(x)| + \frac{1}{2}\rho(x,t).
    \end{align*}
    Thus, by (\ref{eq: compare lambda}) we find
    \begin{align*}
        \rho(x',t')^{2} = \lambda(t')^{2} + |F_{t'}(x')|^{2}\leq K^{2}\lambda(t)^{2} + 2|F_t(x)|^{2} + \frac{1}{2}\rho(x,t)^{2}\leq \left(\frac{5}{2} + K^{2}\right)\rho(x,t)^{2}.
    \end{align*}
    On the other hand, if $|F_t(x)|\leq\lambda(t)$,
    \begin{align*}
        \rho(x',t')^{2}\geq\lambda(t')^{2}\geq \frac1{K^2}\lambda(t)^{2}\geq\frac{1}{2K^{2}}\rho(x,t)^{2},
    \end{align*}
    and if $|F_t(x)|>\lambda(t)$, then $|F_t(x)| - \frac{1}{2}\rho(x,t)>\left(1 - \frac{1}{\sqrt{2}}\right)|F_t(x)|$. Thus,
    \begin{align*}
        \rho(x',t')^{2}\geq|F_{t'}(x')|^{2}\geq \left(|F_t(x)| - \frac{1}{2}\rho(x,t)\right)^{2}\geq \left(1 - \frac{1}{\sqrt{2}}\right)^{2}|F_t(x)|^{2}>\frac{1}{2}\left(1 - \frac{1}{\sqrt{2}}\right)^{2}\rho(x,t)^{2}.
    \end{align*}
    Hence, (\ref{eq: compare rho}) is proved.

    Finally, since $R(t)$ is non-increasing, we have
    \begin{multline*}
        |F_{t'}(x')|\leq |F_t(x)| + \frac{1}{2}\rho(x,t) \leq C^{-1}R(t) + \frac{1}{2}\sqrt{\lambda(t)^{2} + C^{-2}R(t)^{2}} \\
        \leq\left(C^{-1} + \frac{1}{2}\sqrt{\frac{\lambda(t)^{2}}{R(t)^{2}}+ C^{-2}} \right)R(t').
    \end{multline*}
    By choosing $t_{0}$ close enough to $T$ we can make $C^{-1} + \frac{1}{2}\sqrt{\frac{\lambda(t)^{2}}{R(t)^{2}}+ C^{-2}} \leq 1$. This proves (\ref{eq: compare global radius}).
\end{proof}

The main pointwise mean curvature bound Thereom \ref{cor: pointwise mean curvature bound} now follows directly from combining Theorem \ref{thm angle est}, Theorem \ref{thm deriv angle est}, and Lemma \ref{lem: compare weights}.

\begin{proof}[Proof of Theorem~\ref{cor: pointwise mean curvature bound}]
    By Theorem~\ref{thm angle est}, for $(x, t) \in L \times [0, T)$ with $|F_t(x)|\leq R(t)$, we have
    \begin{align*}
        |\theta(x,t) - \overline{\theta}|\leq CM\rho(x,t)^{-\alpha}\lambda(t)^{\beta}
    \end{align*}
    where
    \begin{equation*}
        M:=
        \sup_{t' \in [0, T)} \sup_{\substack{x' \in L, \\ C^{-1} R(t') \le | F_{t'}(x') | \le  R(t')}} \frac{\rho(x',t')^\alpha}{\lambda(t')^\beta } \cdot  | \theta(x', t') - \overline \theta|.
    \end{equation*}
    It follows from Lemma~\ref{lem: compare weights} that if $(x,t)\in L \times [t_{0}, T)$ with $|F_t(x)|\leq C^{-1}R(t)$, then
    \begin{align*}
        \sup_{P_{L}\left(x, t, \frac{1}{2}\rho(x,t)\right)}|\theta - \overline{\theta}|\leq C M \rho(x, t)^{-\alpha}\lambda(t)^{\beta}.
    \end{align*}
    The result now follows from Theorem~\ref{thm deriv angle est}.
\end{proof}

\section{Applications}

In this section, we apply the mean curvature bound from Theorem \ref{cor: pointwise mean curvature bound} to deduce the existence of both finite-time and infinite-time LMCF singularities with uniformly bounded mean curvature (Corollary \ref{main cor intro}).

\subsection{Finite-Time Singularities}

\cite{SS26II} constructed a collection of Lagrangian mean curvature flows that form finite-time singularities with precise curvature blow-up rates, descriptions of the singularity models, and estimates on the difference from the singularity models.
In this subsection, we prove a subcollection of these solutions have uniformly bounded mean curvature (Corollary \ref{main cor finite-time sings sect}).
Let us first recall the main results and estimates from \cite{SS26II}.

\begin{theorem}[\cite{SS26II}] \label{thm SS26II}
    Let $n, K \in \mathbb N$ be such that 
        $$(n \ge 4 \text{ and } K \ge 2) \qquad \text{or} \qquad ( n = 3 \text{ and } K \in \{ 2,3,4\}).$$
    Then there exists a smooth, non-compact, properly embedded, almost calibrated, $O(n)$-invariant LMCF $L_K(t) \subset \C^n$ defined for $0 \le t < T = T(n, K) < \infty$ 
    such that 
    \begin{enumerate}
        \item \label{thm SS26II, item singularity at origin}
        the flow $L_K(t)$ forms a singularity at the origin $\mathbf 0 \in \C^n$ at time $T < \infty$,
        and
        \item \label{thm SS26II, item curv blow-up rate}
        the second fundamental form blows up as $\sup_{L_K(t)} |\mathbf A_{L_K(t)}| \sim \frac1{(T-t)^{K/2}}$, that is,
            $$0 < \liminf_{t \nearrow T} \sup_{\mathbf x \in L_K(t)} ( T-t)^{K/2} | \mathbf A_{L_K(t)}| 
            \le \limsup_{t \nearrow T} \sup_{\mathbf x \in L_K(t)} ( T-t)^{K/2} | \mathbf A_{L_K(t)}|  < \infty.$$        
    \end{enumerate}

    The $O(n)$-invariant rescaled LMCF $e^{\tau /2} L_K(T - e^{-\tau})$ is given by the graph of a profile function $h : \R \times [-\ln T, \infty) \to \R$, $h = h(s, \tau)$, that is even in $s$ for every $\tau$ and takes the form
        $$h(s, \tau) = f_{a(\tau)}(s) + \partial_s u(s, \tau), \quad 
        u(s, \tau) = \sum_{i=1}^K b_i (\tau) \left( \phi_{i,a(\tau)}(s) - \phi_{0, a(\tau)}(s) \right) + v(s, \tau)$$
    where
    \begin{enumerate}
        \item $f_a(s) = a f_1 (s/a)$ is the profile function of the $O(n)$-invariant Lawlor neck with $f_a(0) = a > 0$
        and

        \item $\phi_{i,a(\tau)}$ $(0 \le i \le K)$ denote the eigenfunctions given in \cite{SS26I}*{Theorem \ref{I-main thm intro}} 
    \end{enumerate}
    such that the terms satisfy the following estimates for all $\tau \in [-\ln T, \infty)$
    \begin{enumerate}
        \item ($a(\tau)$ Estimates)
        $$\frac12 e^{\frac{1-K}2 \tau} \le a(\tau) \le 2 e^{\frac{1-K}2 \tau} < 1 \quad \text{and} \quad \left| \frac{da}{d\tau} \right| \le (K-1) \cdot a(\tau),$$

        \item ($b_i(\tau)$ Estimates)
            $$\sum_{i=1}^K |b_i(\tau) | \le 2K \cdot a(\tau)^2$$ 
        and
        \item (Weighted $C^2$ Estimates)
        \begin{align}
            \label{thm SS26II, eqn eigfunc inner region est}
            \max_{1 \le i \le K }   \| \phi_{i,a(\tau)} - \phi_{0,a(\tau)} \|_{C^2_{\left(1 + \frac \gamma 2\right) } (\{ |s|< s_0 \})} &\lesssim_{n, K, \gamma, s_0} 1,  \\ 
            \| v( \cdot, \tau) \|_{C^2_{(\gamma)} \left( \{ |s| < s_0 \} \right) } &\le e^{-\kappa_{in} \tau}, \text{ and} \\
            \| u ( \cdot, \tau) \|_{C^2_{(0)} \left( \{ \frac{s_0}2 < |s| < 2\Upsilon \} \right) }  &\lesssim_{n, K, s_0, \Upsilon} e^{-( K- 1) \tau} 
        \end{align}
        where $\| \cdot \|_{C^2_{(\omega)}(\Omega)}$ denotes the ($\tau$-dependent) weighted $C^2$ norm defined by
            $$\| u(\cdot, \tau) \|_{C^2_{(\omega)} (\Omega) } \doteqdot \sum_{j=0}^2 \sup_{s \in \Omega} \tilde \rho(s, \tau)^{-\omega +j} |\partial_s^j u|(s, \tau), \qquad \tilde \rho(s, \tau) \doteqdot \sqrt{a(\tau)^2 + s^2} ,$$
        and $\gamma , \kappa_{in}, \Upsilon, s_0$ are parameters that satisfy
            $$-1 < \gamma < 0, \qquad \kappa_{in} = ( 2+ \tilde \kappa_{in} ) \frac{K-1}2 \text{ with }  - \gamma \le \tilde \kappa_{in}, \qquad 
            \text{and} \qquad 0 < s_0 < 1 < \Upsilon < \infty.$$
        \end{enumerate}
        Additionally, the profile function $\tilde h(x,t) = \sqrt{T-t} \cdot  h ( x/ \sqrt{T-t}, -\ln(T-t))$ of the $O(n)$-invariant LMCF $L_K(t)$ satisfies the following weighted $C^2$-estimates with the parameter $\kappa_{out} = \frac78 (K-1)$: for any $0 < \epsilon < 1$,
        \begin{equation} \label{thm SS26II, eqn profile est outer region}
            \sup_{t \in [0, T)} \sup_{|x| > \epsilon \sqrt{T-t}} \sum_{j=0}^2 \left( |x|^{-1+j} + |x|^{-2 \kappa_{out} -1 + j} \right) \left|\partial_x^j \left( \tilde h (x, t) - \cot \left(\frac{\pi}{2n} \right) |x| \right) \right| < \infty .
        \end{equation}
\end{theorem}

\begin{remark}
    Theorem \ref{thm SS26II} differs from the statement of \cite{SS26II}*{Theorem \ref{II-main thm paper II intro}} in that Theorem \ref{thm SS26II} contains additional estimates for the profile function of the solution.
    These estimates are proven in \cite{SS26II}*{Section \ref{II-Section Box Argument}}.
    Indeed, the proof \cite{SS26II}*{Theorem \ref{II-main thm paper II intro}} follows by appealing to \cite{SS26II}*{Corollary \ref{II-cor main box cor}}, which implies all the estimates in Theorem \ref{thm SS26II} above except \eqref{thm SS26II, eqn eigfunc inner region est}.
    Specifically, the estimates for $a, b_i, u, v$ follow from the definition of $\mathcal B$ \cite{SS26II}*{Definition \ref{II-defn box}} and estimate \eqref{thm SS26II, eqn profile est outer region} follows from invoking \cite{SS26II}*{Lemmas \ref{II-lem C^2 Profile Ests on outer-parab overlap}--\ref{II-lem C^2 Profile Ests in outer region}} in the proof of \cite{SS26II}*{Theorem \ref{II-main thm paper II intro}}.
    The remaining estimate \eqref{thm SS26II, eqn eigfunc inner region est} is obtained in \cite{SS26II}*{Equation \eqref{II-proof preserving inner Holder condition, proof claim 4, eqn 5}}.
\end{remark}

These estimates from \cite{SS26II} already suffice to prove the mean curvature stays bounded outside parabolic scales $|\mathbf x | \gtrsim \sqrt{T-t}$.

\begin{lem} \label{lem H bounds outside parab region}
    Let $n, K \in \mathbb N$ be such that
        $$(n \ge 4 \text{ and } K \ge 2) \qquad \text{or} \qquad ( n = 3 \text{ and } K \in \{ 2,3,4\}).$$
    Let $L_K(t) \subset \C^n$ ($t \in [0, T)$) be the $O(n)$-invariant LMCF solution from Theorem \ref{thm SS26II}.

    Then, for any $0 < \epsilon < 1$, the mean curvature $\mathbf H = \mathbf H_{L_K(t)}$ satisfies
    \begin{equation}
        \sup_{t \in [0, T)} \sup_{\substack{\mathbf x \in L_K(t),\\|\mathbf x| \ge \epsilon \sqrt{T-t}}} | \mathbf H | < \infty.
    \end{equation}
\end{lem}
\begin{proof}
    Let $\tilde h (x,t) = \sqrt{T - t} \cdot h ( x/ \sqrt{T - t} , - \ln (T-t))$ denote the profile function for the $O(n)$-invariant LMCF $L_K(t)$.
    Then the magnitude of the mean curvature at the point corresponding to $(x,t)$ can be expressed as
    \begin{equation*}
        | \mathbf H| = \frac1{\sqrt{1 + \tilde h_x^2}} \left| \frac{\tilde h_{xx}}{1 + \tilde h_x^2} + ( n-1) \frac{ x \tilde h_x - \tilde h}{x^2 + \tilde h^2} \right| 
    \end{equation*}
    (see e.g. \cite{SS26I}*{Equation \eqref{I-eqn defn H[h]}}).
    Writing $\overline c_0 := \cot \left( \frac \pi{2n}\right) $ and $\hat h(x,t) := \tilde h(x,t) - \overline c_0 |x| $ to be the difference from the cone profile, it follows that
    \begin{equation} \label{proof H bound outside parab region, eqn 1}
        | \mathbf H| = \frac1{\sqrt{1 + (\overline c_0 +\hat h_x)^2}} \left| \frac{\hat h_{xx}}{1 + (\overline c_0 +\hat h_x)^2} + ( n-1) \frac{ x \hat h_x - \hat h}{x^2 + (\overline c_0 x + \hat h)^2} \right| .
    \end{equation}
    on the region where $x>0$.
    
    With \eqref{proof H bound outside parab region, eqn 1}, the bounds for $\hat h$ from \eqref{thm SS26II, eqn profile est outer region} imply
    \begin{align*}
        \sup_{t \in [0, T)} \sup_{|x| \ge \epsilon \sqrt{T-t}} |\mathbf H|  
        &\lesssim \sup_{t \in [0, T)} \sup_{ \epsilon \sqrt{T-t}\le |x| } \left( |\hat h_{xx} | + \left| \frac1x \hat h_x \right| + \left| \frac1{x^2} \hat h \right|  \right) \\
        &\lesssim 
        \sup_{t \in [0, T)} \sup_{ \epsilon \sqrt{T-t}\le |x| \le \epsilon} |x|^{2 \kappa_{out}-1} 
        + \sup_{t \in [0, T)} \sup_{ \epsilon \le |x|} |x|^{-1} 
    \end{align*}
    for any $0 < \epsilon < 1$,
    where the implicit constants depend on $n, K, \epsilon, T$ and the finite supremum in \eqref{thm SS26II, eqn profile est outer region}.
    Since $\kappa_{out} = \frac78 (K-1) \ge \frac12$, the last line is finite and thus
        $$\sup_{t \in [0, T)} \sup_{x \ge \epsilon \sqrt{T-t}} |\mathbf H| < \infty \qquad \text{for any } 0 < \epsilon < 1.$$
    Finally, the $O(n)$-invariance of the flow and the estimates for $\tilde h$ \eqref{thm SS26II, eqn profile est outer region}
    imply the same finiteness among all points $\mathbf x \in L_K(t)$ with $|\mathbf x| \ge \epsilon \sqrt{T-t}$.
\end{proof}

It remains to bound the mean curvature within parabolic scales $|\mathbf x| \lesssim \sqrt{T-t}$. 
To do so, we will appeal to the general results developed in Section \ref{Section Ests at Small Scales}.
The following proposition first verifies Assumption \ref{blow-up assumptions} holds in this case.
The proof is technical but essentially follows from the estimates listed in Theorem \ref{thm SS26II} above.

\begin{proposition} \label{prop assumptions hold finite time}
    Let $n, K \in \mathbb N$ be such that 
        $$(n \ge 4 \text{ and } K \ge 2) \qquad \text{or} \qquad ( n = 3 \text{ and } K \in \{ 2,3,4\}).$$
    Let $L_K(t) \subset \C^n$ be the LMCF solution from Theorem \ref{thm SS26II}.  
    If $R(t) = \sqrt{T-t}$ and $\lambda(t) = (T - t)^{K/2}$,
    then $(L_K(t), R(t), \lambda(t))$ satisfies Assumption \ref{blow-up assumptions}.
\end{proposition}
\begin{proof}    
    Condition \eqref{scale assumptions} of Assumption \ref{blow-up assumptions} follows immediately from the choice of $R(t), \lambda(t)$.

    Condition \eqref{assumption, blow-up to Lawlor neck} of Assumption \ref{blow-up assumptions} follows from applying \cite{Wood24}*{Theorem 1.2}.


    For Condition \eqref{assumption, blow-up to cone}, consider sequences $t_i \to T$ and $\mu_i$ such that $(T-t_i)^{K/2} \ll \mu_i \ll \sqrt{T - t_i}$, and 
    consider the sequence of rescaled flows
        $$\hat L_i(\hat \tau) := \mu_i^{-1}  L_K\left(t_i + \hat \tau \cdot \mu_i^2 \right). $$
    Since $e^{\tau/2} L_K ( T - e^{-\tau})$ is given by the graph of the profile function $h(s, \tau) = f_{a(\tau)} (s) + \partial_s u(s, \tau)$ as in Theorem \ref{thm SS26II}, it follows that $\hat L_i(\hat \tau)$ is given by the graph of the profile function 
    \begin{align*}
        \hat h_i (\xi, \hat \tau) 
        ={}& \mu_i^{-1} \sqrt{ T - t_i - \hat \tau \mu_i^2 } \cdot \, h \left( \frac{\xi}{\mu_i^{-1} \sqrt{ T - t_i - \hat \tau \mu_i^2 }}, - \ln ( T - t_i - \hat \tau \mu_i^2 ) \right) \\
        ={}& f_{\hat a_i(\hat \tau) } (\xi) 
        +\partial_\xi \hat u_i (\xi, \hat \tau)
    \end{align*}
    where
    \begin{gather}
        \label{proof assumptions hold finite time, eqn 1}
        \hat a_i(\hat \tau) := a (-\ln ( T - t_i - \hat \tau \mu_i^2)) \cdot \mu_i^{-1} \sqrt{ T - t_i - \hat \tau \mu_i^2 } \qquad \text{ and}\\
        \label{proof assumptions hold finite time, eqn 2}
        \hat u_i (\xi, \hat \tau) :=  \mu_i^{-2}  ( T - t_i - \hat \tau \mu_i^2 ) \cdot u \left( \frac{\xi}{ \mu_i^{-1} \sqrt{ T - t_i - \hat \tau \mu_i^2} } , - \ln ( T - t_i - \hat \tau \mu_i^2) \right).
    \end{gather}
    
    Because $a(\tau) \sim e^{\frac{1-K}2 \tau}$ (Theorem \ref{thm SS26II}), we compute from \eqref{proof assumptions hold finite time, eqn 1} that
    \begin{align*}
        \hat a_i (\hat \tau) &\sim ( T - t_i - \hat \tau \mu_i^2)^{\frac K2} \mu_i^{-1}. 
    \end{align*}
    Using that $K \ge 2$, $t_i \nearrow T$, and $(T -t_i)^{\frac K2} \ll \mu_i \ll \sqrt{ T - t_i}$, it follows that $\hat a_i (\hat \tau) \xrightarrow[i \to \infty]{C^0_{\loc}(\R)} 0 $.
    Using also $f_{\hat a_i (\hat \tau)} (\xi) = \hat a_i(\hat \tau) f_1 ( \xi / \hat a_i( \hat \tau)) $ and the asymptotics of the Lawlor neck profile function $f_1$ (see e.g. \cite{SS26I}*{Lemma~\ref{I-lem: profile function of Lawlor neck}}), we therefore deduce
    \begin{equation} \label{proof assumptions hold finite time, eqn 3}
        \sup_{-A^2 < \hat \tau \le 0} \left\| f_{\hat a_i (\hat \tau) }(\xi)  - \cot \left( \frac \pi{2n} \right) \cdot  |\xi| \right\|_{C^1( \{ A^{-1} < |\xi| < A\} ) } \xrightarrow[i \to \infty]{} 0
    \end{equation}
    for any constant $1 < A < \infty$.
    Note that $\cot \left( \frac \pi {2n} \right) |\xi|$ is the profile function for an $O(n)$-invariant special Lagrangian transverse pair of planes.

\begin{claim} \label{proof assumptions hold finite time, claim 1}
        For any $1 < A < \infty$, 
            $$\sup_{-A^2 < \hat \tau \le 0 } \| \hat u_i (\xi, \hat \tau) \|_{C^2 ( \{ A^{-1} < |\xi| < A \} )}  \xrightarrow[i \to \infty]{}0.$$
    \end{claim}
    \begin{claimproof}
        Fix $1 < A < \infty$.
        Because $\mu_i \ll \sqrt{T - t_i}$,  
            $$\mu_i^{-1} \sqrt{ T- t_i - \hat \tau \mu_i^2} \xrightarrow[i \to \infty]{} + \infty \qquad \text{uniformly in } \hat \tau \in (-A^2, 0].$$
        Hence, when $i \gg 1$ is sufficiently large,
        \begin{equation} \label{proof assumptions hold finite time, proof claim 1, eqn 1}
            \frac{ |\xi|}{\mu_i^{-1} \sqrt{ T- t_i - \hat \tau \mu_i^2}} < s_0 
            \qquad \text{for all } A^{-1} < |\xi| < A \text{ and } -A^2 < \hat \tau \le 0,
        \end{equation}
        where $s_0>0$ is as in Theorem \ref{thm SS26II}.

        Using $( T - t_i )^{K/2} \ll \mu_i$ and $a(\tau) \sim e^{\frac{1-K}2 \tau}$ (Theorem \ref{thm SS26II}), one can obtain 
        \begin{multline} \label{proof assumptions hold finite time, proof claim 1, eqn 2}
            \frac{|\xi|}{\mu_i^{-1} \sqrt{ T - t_i - \hat \tau \mu_i^2}}
            \le \sqrt{ a^2 \left( - \ln ( T - t_i - \hat \tau \mu_i^2) \right) + \frac{|\xi|^2}{\mu_i^{-2} ( T - t_i - \hat \tau \mu_i^2)}} \\
            \le  2  \frac{ |\xi|}{\mu_i^{-1} \sqrt{ T - t_i - \hat \tau \mu_i^2}}
        \end{multline}
        for all $A^{-1} < |\xi| < A$ and $-A^2 < \tau \le 0$, when $i \gg 1$ is sufficiently large.   

        By the definition of $\hat u_i$ \eqref{proof assumptions hold finite time, eqn 2} and $u$ (Theorem \ref{thm SS26II}),
        \begin{equation} \label{proof assumptions hold finite time, proof claim 1, eqn 3}
            \partial_\xi^j \hat u_i (\xi, \hat \tau) 
            = \mu_i^{j-2} ( T - t_i - \hat \tau \mu_i^2 )^{1 - \frac j2} \left( \sum_{i=1}^K b_i(\tau) \partial_s^j \left(\phi_{i,a(\tau)} - \phi_{0, a(\tau)}\right )  + \partial_s^j v \right) 
        \end{equation}
        where the functions on the right-hand side are evaluated at $s = \frac \xi {\mu_i^{-1} \sqrt{ T - t_i - \hat \tau \mu_i^2} }$ and time $ \tau = -\ln ( T - t_i - \hat \tau \mu_i^2)$. 
        We estimate the contributions of these terms for $0 \le j \le 2$ as $i \to \infty$.
        In what follows, all estimates are uniform in $|\xi| \in (A^{-1} , A)$ and $\hat \tau \in (-A^2, 0]$.
        The contributions of the $\phi_{i,a}$ terms can be estimated as
        {\allowdisplaybreaks 
        \begin{align*}
            & \qquad \mu_i^{j-2} ( T - t_i - \hat \tau \mu_i^2 )^{1 - \frac j2} \sum_{i=1}^K |b_i ( \tau)| | \partial_s^j ( \phi_{i,a(\tau)} - \phi_{0,a(\tau)} ) | \\
            &\lesssim_{K} \mu_i^{j-2} ( T - t_i - \hat \tau \mu_i^2 )^{1 - \frac j2 + K - 1}  \sum_{i=1}^K  | \partial_s^j ( \phi_{i,a} - \phi_{0,a} ) | \\
            & \qquad \qquad  \left( \text{since } |b_i| \lesssim_K a^2 \sim e^{(1-K)  \tau} \right) \\
            &\le \mu_i^{j-2} ( T - t_i - \hat \tau \mu_i^2 )^{K - \frac j2 } 
            \tilde \rho(s,\tau)^{1+ \frac \gamma 2-j} \sum_{i=1}^K \| \phi_{i,a} - \phi_{0,a} \|_{C^2_{\left(1+ \frac \gamma 2\right)}( \{ |s| < s_0 \} ) }
            && (0 \le j \le 2, \, \eqref{proof assumptions hold finite time, proof claim 1, eqn 1}) \\
            &\lesssim_{n, K, s_0, \gamma} \mu_i^{j-2} ( T - t_i - \hat \tau \mu_i^2 )^{K - \frac j2 } 
            \tilde \rho(s,\tau)^{1+ \frac \gamma 2 -j}  
            && ( \text{by \eqref{thm SS26II, eqn eigfunc inner region est}})\\
            &\lesssim_{A, \gamma} \mu_i^{j-2 +1+ \frac \gamma 2-j} ( T - t_i - \hat \tau \mu_i^2 )^{K - \frac j2 - \frac12 -\frac\gamma 4+ \frac j2 } 
            && ( \text{by \eqref{proof assumptions hold finite time, proof claim 1, eqn 2}}) \\
            &= \mu_i^{-1+\frac \gamma2} ( T - t_i - \hat \tau \mu_i^2 )^{K  - \frac12 -\frac \gamma 4  } \\
            &\ll ( T - t_i)^{ \frac K2(-1+ \frac \gamma 2)} ( T - t_i  )^{K  - \frac12  - \frac \gamma 4} \\
            & \qquad \qquad \left( \text{since } ( T - t_i)^{K/2} \ll \mu_i \ll \sqrt{ T - t_i} , \, \gamma < 0, \, \text{and } |\tau| < A^2 \right) \\
            &= ( T - t_i )^{\frac K2 - \frac 12+ \frac \gamma 4 (K-1)} \\
            &\ll ( T - t_i )^{\frac K2 - \frac12 - \frac14 (K-1)}
            =( T - t_i )^{\frac K4 - \frac14 } \ll 1 \\
            & \qquad  \qquad ( \text{since } K \ge 2, \, t_i \nearrow T, \, -1 < \gamma < 0) .
        \end{align*}}
        The contributions of the $v$ terms can be estimated as
        {\allowdisplaybreaks
        \begin{align*}
            & \quad \mu_i^{j-2} ( T - t_i - \hat \tau \mu_i^2)^{1 - \frac j2} | \partial_s^j v| \\
            &\le \mu_i^{j-2} ( T - t_i - \hat \tau \mu_i^2)^{1 - \frac j2} \tilde \rho(s,\tau)^{\gamma - j} \|  v\|_{C^2_{(\gamma)} ( \{ |s | < s_0 \} ) } 
            && ( 0 \le j \le 2, \, \eqref{proof assumptions hold finite time, proof claim 1, eqn 1}) \\
            &\le \mu_i^{j-2} ( T - t_i - \hat \tau \mu_i^2)^{1 - \frac j2} \tilde \rho(s,  \tau)^{\gamma - j} e^{- \kappa_{in}  \tau}
            && ( \text{by Theorem \ref{thm SS26II}} )\\
            &\lesssim_{A, \gamma} \mu_i^{-2+ \gamma} ( T - t_i - \hat \tau \mu_i^2)^{1 - \frac \gamma 2}  e^{- \kappa_{in}  \tau}
            && ( \text{by } \eqref{proof assumptions hold finite time, proof claim 1, eqn 2} ) \\
            &= \mu_i^{-2+ \gamma} ( T - t_i - \hat \tau \mu_i^2)^{1 - \frac \gamma 2+ \kappa_{in}} 
            && \left(  \tau  = - \ln \left( T - t_i - \hat \tau \mu_i^2\right) \right)\\
            &\ll ( T - t_i)^{\frac K2 ( - 2 + \gamma) } ( T - t_i )^{1 - \frac \gamma 2+ \kappa_{in}} = 
            ( T - t_i)^{\frac K2 ( - 2 + \gamma)+  1 - \frac \gamma 2+ \kappa_{in} }\\
            &\qquad  \quad \left( \text{since }  ( T - t_i )^{\frac K2 } \ll \mu_i \ll \sqrt{T - t_i} , \, \gamma <  0, \, \text{and } |\tau| \le A^2 \right) .
        \end{align*}}
        Writing $\kappa_{in} = ( 2  +\tilde \kappa_{in} ) \frac{K-1}2$ as in Theorem \ref{thm SS26II}, the exponent in the last line is 
        $$\frac{K}2 ( - 2 + \gamma) + 1 - \frac \gamma 2 + \kappa_{in}
        = \frac{K-1}2 ( \gamma + \tilde \kappa_{in} ) \ge 0$$
        since $K \ge 2$ and $\tilde \kappa_{in} \ge - \gamma$.
        Thus, 
        \begin{equation}
            \mu_i^{j-2} ( T - t_i - \hat \tau \mu_i^2)^{1 - \frac j2} | \partial_s^j v| 
            \ll  ( T - t_i)^{\frac K2 ( - 2 + \gamma)+  1 - \frac \gamma 2+ \kappa_{in} }
            = ( T - t_i)^{ \frac{K-1}2 ( \gamma + \tilde \kappa_{in} ) } 
            \le 1.
        \end{equation}
        Combining these estimates for the contributions of the two terms in $\partial_\xi^j \hat u_i$ \eqref{proof assumptions hold finite time, proof claim 1, eqn 3} thereby proves the claim.        
    \end{claimproof}

    Claim \ref{proof assumptions hold finite time, claim 1} and \eqref{proof assumptions hold finite time, eqn 3} show that the sequence $\hat L_i(\hat \tau)$ of mean curvature flows converges to a special Lagrangian transverse pair of planes, and the convergence is $C^0_{\loc}$ in $\hat \tau \in (-\infty, 0]$ and $C^1_{\loc}$ in $\xi \in \C^n \setminus \{ \mathbf 0\}$.
    This convergence is strong enough to apply White's Local Regularity Theorem \cite{White05} to improve the convergence to $C^\infty_{\loc}( \C^n \setminus \{ \mathbf 0 \} \times (-\infty, 0])$-convergence to the pair of planes.
    This proves that Condition \ref{assumption, blow-up to cone} of Assumption \ref{blow-up assumptions} holds, and thereby completes the proof that $(L_K(t), R(t), \lambda(t))$ satisfies Assumption \ref{blow-up assumptions}.
\end{proof}

\begin{proposition} \label{prop finite-time sing, angle decay}
    Let $n, K \in \mathbb N$ be such that 
        $$(n \ge 4 \text{ and } K \ge 2) \qquad \text{or} \qquad ( n = 3 \text{ and } K \in \{ 2,3,4\}).$$
    Let $L_K(t) \subset \C^n$ be the LMCF solution from Theorem \ref{thm SS26II}.  

    Then, for any $0 < \epsilon < 1$, the Lagrangian angle $\theta = \theta_{L_K(t)}$ satisfies
    \begin{equation}
        \sup_{t \in [0,T)} \sup_{\substack{\mathbf x \in L_K(t),\\ \epsilon \sqrt{T-t}\le |\mathbf x|\le \sqrt{T-t}}} \frac{|\theta - (n-1)\frac \pi 2 |}{(T-t)^{K-1}} < \infty.
    \end{equation}
\end{proposition}

\begin{proof}
    By Theorem \ref{thm SS26II}, the rescaled flow $e^{\tau/ 2} L_K ( T - e^{-\tau})$ is given by the graph of a profile function $h = f_{a(\tau)}(s) + \partial_s u (s, \tau)$.
    The Lagrangian angle $\theta = \theta_{L(T - e^{-\tau})}$ at $\mathbf x = e^{-\tau/2} ( s + i h(s, \tau), 0, \dots, 0) \in L_K ( T - e^{-\tau} ) \subset \C^n$
    can then be expressed in terms of the profile function $h$ as
        $$\theta = \arctan (h_s) + (n-1) \arctan \left( \frac h s \right) $$
    (see e.g. \cite{SS26I}*{Lemma \ref{I-lem: potentials of H and X}}).
    For $h = f_a + \partial_s u$, \cite{SS26I}*{Proposition \ref{I-prop: PDE for the potential}} yields
    \begin{gather}
        \label{prop finite-time sing, angle decay, eqn 1}
        \theta - (n-1) \frac \pi 2 = L_{s, a} u + Q_{s,a}(u)  , \qquad  a = a(\tau) \\
        \nonumber \text{where } L_{s,a} u= \frac{u_{ss}}{1 + (\partial_s f_{a})^2} + (n-1) \frac{ s u_s}{s^2 + f_{a}^2}. 
    \end{gather}
    
    \begin{claim} \label{proof finit-time sing, angle decay, claim 1}
        For any $0 < \epsilon < 1$,
        \begin{equation}
            \sup_{\tau \ge -\ln T} \sup_{\epsilon \le |s| \le 2 \Upsilon} e^{(K-1)\tau} |L_{s,a} u | < \infty ,
        \end{equation}        
        where $\Upsilon > 1$ is given in Theorem \ref{thm SS26II}.
    \end{claim}
    \begin{claimproof}
        Since $a(\tau) \le 1$ for all $\tau$ and the functions $f_a (s)$ smoothly converge to the profile function $\cot\left( \frac \pi {2n} \right) |s|$ of the special Lagrangian pair of planes as $a \to 0$ in the region $\epsilon \le |s| \le 2 \Upsilon$, 
        there exists $C = C(n, \epsilon , \Upsilon) > 1$ such that 
        \begin{equation*}
            C^{-1} < | \partial_s f_a | , \left| \frac {f_a} s \right| < C \qquad \forall 0 < a \le 1, \, \epsilon \le |s| \le 2 \Upsilon.
        \end{equation*}
        Hence,
        \begin{equation} \label{proof finit-time sing, angle decay, claim 1, eqn 0}
            | L_{s,a} u | \lesssim_{n, \epsilon, \Upsilon} |u_{ss} | + |u_s| \qquad \forall \tau \ge - \ln T, \, \epsilon \le |s| \le 2 \Upsilon.
        \end{equation}
        
        Note also that $\epsilon \le \tilde \rho(s,\tau) = \sqrt{ a^2( \tau) + s^2} \le \sqrt{1 + 4 \Upsilon^2} $ in this region, and so the weighted H{\"o}lder norms from Theorem \ref{thm SS26II} are equivalent to standard H{\"o}lder norms.
        The H{\"o}lder estimates from Theorem \ref{thm SS26II} therefore imply, when $\max \{ \epsilon, \frac{s_0}2 \} \le |s| \le 2 \Upsilon$,
        \begin{equation} \label{proof finit-time sing, angle decay, claim 1, eqn 1}
            |L_{s,a} u | \lesssim_{n, \epsilon , \Upsilon} |u_{ss} | + |u_s | \lesssim_{n, K, s_0, \Upsilon}
             e^{-(K-1) \tau }
        \end{equation}
        Similarly, when $\epsilon \le |s| \le \frac{s_0}2$, the estimates from Theorem \ref{thm SS26II} imply
        \begin{gather} \label{proof finit-time sing, angle decay, claim 1, eqn 2} \begin{aligned}
            & \quad |L_{s,a} u | \\
            &\lesssim_{n, \epsilon, \Upsilon, \gamma} \sum_{i=1}^K |b_i| \| \phi_{i,a} - \phi_{0,a} \|_{C^2_{(1 + \frac \gamma 2 )}(\{ |s| < s_0\})}
            + \| v \|_{C^2_{(\gamma)}(\{ |s| < s_0\})} 
            && ( \text{by \eqref{proof finit-time sing, angle decay, claim 1, eqn 0}}) \\
            &\lesssim_{n,K,\gamma, s_0} a^2 + e^{-\kappa_{in} \tau} 
            && (\text{by Theorem \ref{thm SS26II}})\\
            &\lesssim_{\kappa_{in}, K, T} e^{- (K-1) \tau}
        \end{aligned} \end{gather}
        where the last inequality uses the facts $a \sim e^{-\frac{K-1}2 \tau}$ and $\kappa_{in} > K-1$ (Theorem \ref{thm SS26II}).
        Combining estimates \eqref{proof finit-time sing, angle decay, claim 1, eqn 1} and \eqref{proof finit-time sing, angle decay, claim 1, eqn 2} proves the claim.
    \end{claimproof}

    \begin{claim} \label{proof finit-time sing, angle decay, claim 2}
        For any $0 < \epsilon < 1$,
        \begin{equation}
            \sup_{\tau \ge -\ln T} \sup_{\epsilon \le |s| \le 2 \Upsilon} e^{(K-1)\tau} |Q_{s,a} (u) | < \infty 
        \end{equation}
        where $\Upsilon > 1$ is given in Theorem \ref{thm SS26II}.
    \end{claim}
    \begin{claimproof}
        \cite{SS26II}*{Lemma \ref{II-lem: error term C0}} obtains that
        \begin{equation*}
            |Q_{s,a}(u) | \le n C \left( u_{ss}^2 + \frac{u_s^2}{s^2} \right)
            \qquad \forall s \in \R, \, a > 0
        \end{equation*}
        for some universal constant $C>0$.
        The claim then follows by similar logic as in the proof of Claim \ref{proof finit-time sing, angle decay, claim 1} using the $C^2$-estimates for $u$ from Theorem \ref{thm SS26II}.
    \end{claimproof}

    Claims \ref{proof finit-time sing, angle decay, claim 1} and \ref{proof finit-time sing, angle decay, claim 2} and equation \eqref{prop finite-time sing, angle decay, eqn 1} show that, for any $\epsilon \in (0,1)$,
        $$\sup_{\tau \ge -\ln T } \sup_{\substack{\epsilon \le |s | \le 2 \Upsilon, \\ \mathbf x = e^{-\tau/2} ( s + i h(s, \tau), 0, \dots, 0) \in L_K ( T - e^{-\tau} )}} e^{(K-1)\tau} \left|\theta_{L_K( T - e^{-\tau})} ( \mathbf x) - (n-1) \frac \pi 2 \right|< \infty.$$
    By the $O(n)$-invariance of $L_K$ and the estimates for the profile function $h$ (Theorem \ref{thm SS26II}), these Lagrangian angle bounds on this set of points $(s, \tau)$ imply the Lagrangian angle bounds on the set of points $(\mathbf x,t)$ given in the statement of the proposition.
\end{proof}

We can now apply the results of Section \ref{Section Ests at Small Scales} to deduce the LMCF solutions from \cite{SS26II} with $K \ge 3$ have uniformly bounded mean curvature.
\begin{corollary} \label{main cor finite-time sings sect}
    Let $n, K \in \mathbb N$ be such that 
        $$(n \ge 4 \text{ and } K \ge 3) \qquad \text{or} \qquad ( n = 3 \text{ and } K \in \{ 3,4\}).$$
    Let $L_K(t) \subset \C^n$ be the LMCF solution from Theorem \ref{thm SS26II}.
    Then 
        $$\sup_{t \in [0, T)} \sup_{\mathbf x \in L_K(t)} |\mathbf H| < \infty.$$
\end{corollary}
\begin{proof}
    Lemma \ref{lem H bounds outside parab region} bounds $|\mathbf H|$ on regions where $|\mathbf x | \ge \epsilon \sqrt{T-t}$.
    It therefore suffices to prove
        $$\sup_{t \in [0, T)} \sup_{\substack{\mathbf x \in L_K(t),\\ |\mathbf x| \le C^{-1}\sqrt{T-t}}} |\mathbf H| < \infty \qquad \text{for some constant } C > 1.$$

    Given $n, K \ge 3$ as in the statement, 
    there exists a small positive constant $\alpha \in (0, n-2)$ such that $\alpha ( K - 1) \le K-2$.
    Set $\beta = 1 + \alpha$.

    Proposition \ref{prop assumptions hold finite time} implies Corollary \ref{cor: pointwise mean curvature bound} applies.
    Hence, for some $C> 1$, 
    \begin{align*}
        &\sup_{t \in [0, T)} \sup_{|\mathbf x| \le C^{-1} \sqrt{T-t}} |\mathbf H_{L_K(t)} (\mathbf x)| \\
        \le{}& C  
        +C \frac{ ( T - t)^{\frac K2 \beta  }}{\sqrt{ (T-t)^K+ |\mathbf x|^2}^{1 + \alpha}} 
        \sup_{C^{-1} \sqrt{T - t'} \le |\mathbf x'| \le \sqrt{T - t'}} \frac{\sqrt{ (T-t')^K+ |\mathbf x'|^2}^{\alpha}} {( T - t')^{\frac K2 \beta }}\left|\theta - (n-1) \frac \pi 2\right|
        \\
        \le{}& C 
        + C ( T - t)^{\frac K2 (\beta - 1 - \alpha)   } \sup_{C^{-1} \sqrt{T - t'} \le |\mathbf x'| \le \sqrt{T - t'}} \left( ( T-t')^{ \frac K 2 (  \alpha - \beta )  } + (T - t')^{\frac{ \alpha}{2} - \frac K2 \beta }\right) \left|\theta - (n-1) \frac \pi 2\right|
        \\ 
        ={}& C + C \sup_{C^{-1} \sqrt{T - t'} \le |\mathbf x'| \le \sqrt{T - t'}} \left( (T-t')^{-\frac K2} + (T-t')^{- \frac K2 - \frac \alpha 2 ( K -1)} \right) \left|\theta - (n-1) \frac \pi 2\right|.
    \end{align*}
    By Proposition \ref{prop finite-time sing, angle decay}, this last line will be finite if
        $$\frac K2 \le K-1 \qquad \text{and} \qquad \frac K2 + \frac \alpha2 (K-1) \le K-1.$$
    The former condition holds since $K\ge 2$ and the latter condition holds since $\alpha ( K-1) \le K-2$.
\end{proof}    

\begin{remark}
    Corollary \ref{main cor finite-time sings sect} requires $K \ge 3$.
    It is unclear whether $\mathbf H$ remains bounded in the case $K=2$.
    In this case, Lemma \ref{lem H bounds outside parab region} continues to imply $\mathbf H$ is bounded on spacetime regions where $|\mathbf x| \gtrsim \sqrt{T-t}$, but the proof of Corollary \ref{main cor finite-time sings sect} fails to propagate this estimate to a global bound.
\end{remark}

\subsection{Infinite-Time Singularities}

In \cite{STW}, Su--Tsai--Wood constructed an immortal solution of LMCF
that develops an infinite-time Lawlor neck pinch by a parabolic gluing
argument. We recall the main features of the construction.

\begin{theorem}[\cite{STW}]\label{thm: STW24 solution}
    Let $n\geq 3$, and equip the torus 
    $\mathbb{T}^{2n}=\mathbb{C}^{n}/\Gamma$ with the flat Calabi--Yau
    structure induced by the standard structure on $\mathbb{C}^{n}$.
    Let $X_1$ and $X_2$ be special Lagrangian subtori of phase $0$
    intersecting transversely at a single point $x_*$, and suppose that the union of
    their tangent planes at $x_*$ satisfy the angle criterion such that it can be desingularized by a Lawlor neck $\overline{L}$.

    For all sufficiently small $\varepsilon_0>0$, there is a Lagrangian 
    desingularization $N^{\varepsilon_0}$ of $X:=X_1\cup X_2$, obtained
    by gluing in a Lawlor neck $\overline{L}$ at scale $\varepsilon_0$,
    such that the LMCF $N(t)$ starting from $N^{\varepsilon_0}$ exists
    for all $t\geq 0$. As $t\to\infty$, the flow has the following
    properties:
    \begin{itemize}
        \item \emph{Infinite-time singularity:}
        $N(t)\to X_1\cup X_2$ smoothly away from $x_*$.
        \item \emph{Curvature growth:}
        $\displaystyle\sup_{N(t)}|A_{N(t)}|
        \sim t^{\frac{1}{n-2}}$.
        \item \emph{Singularity model:}
        There are $R>0$ and a positive scale function $\varepsilon(t)$
        satisfying
        \[
            C^{-1}t^{-\frac{1}{n-2}}
            \leq\varepsilon(t)\leq Ct^{-\frac{1}{n-2}}
        \]
        for some $C>1$ and all sufficiently large $t$, such that
        \[
            \varepsilon(t)^{-1}\bigl(N(t)\cap B_R(x_*)\bigr)
            \longrightarrow\overline{L}
        \]
        locally smoothly in the flat coordinates centred at $x_*$.
    \end{itemize}
\end{theorem}

We will show that this flow satisfies
Assumption~\ref{blow-up assumptions} and the weighted angle estimate
required by Corollary~\ref{cor: pointwise mean curvature bound}.
Consequently, its mean curvature remains uniformly bounded despite the
growth of its second fundamental form. We first recall the estimates
from the construction that will be used below.

Choose $r_0>0$ sufficiently small that the flat coordinates centred at
$x_*$ identify $B_{2r_0}(x_*)\subset\mathbb{C}^{n}/\Gamma$ with
$B_{2r_0}(0)\subset\mathbb{C}^{n}$. In these coordinates,
\[
    X\cap B_{2r_0}(0)=\mathcal{C}\cap B_{2r_0}(0),
    \qquad \mathcal{C}:=\Pi_1\cup\Pi_2,
\]
where $\Pi_1$ and $\Pi_2$ are transverse special Lagrangian planes of
phase $0$ satisfying the angle criterion. Let $\overline{L}$ be the
Lawlor neck asymptotic to $\mathcal{C}$ used in the construction.
The approximate solution $N^{\varepsilon(t)}$ is obtained by gluing
$\varepsilon(t)\overline{L}$ into $X$ near the origin. Its relevant
properties are as follows.

\begin{enumerate}
    \item[(a)] \emph{The approximate solution.}
    Outside a compact set, $\overline{L}$ is the Lagrangian graph of
    $\nabla\mathfrak{E}$ over $\mathcal{C}$, with
    \begin{equation}\label{eq: outer Lawlor potential asymptotic}
        |\nabla^k\mathfrak{E}(\mathbf{x})|
        \lesssim |\mathbf{x}|^{2-n-k}
        \qquad\text{as }|\mathbf{x}|\to\infty
    \end{equation}
    for every integer $k\geq 0$.
    Choose a smooth cut-off function
    $\chi:(0,\infty)\to[0,1]$ with $\chi\equiv 0$ on $(0,1]$ and
    $\chi\equiv 1$ on $[2,\infty)$, and set
    \[
        \mathfrak{D}_{\varepsilon}(\mathbf{x})
        :=\bigl[1-\chi(\varepsilon^{-\zeta}|\mathbf{x}|)\bigr]
        \varepsilon^2\mathfrak{E}(\varepsilon^{-1}\mathbf{x}),
        \qquad 0<\zeta<\frac{1}{n+2}.
    \]
    For some fixed $R_L>0$, the approximate solution
    $N^{\varepsilon}$ is represented by the Lagrangian graph of
    $\nabla\mathfrak{D}_{\varepsilon}$ over
    $\mathcal{C}\cap\{R_L\varepsilon<|\mathbf{x}|<r_0\}$.
    In these cone coordinates, $\mathfrak{D}_{\varepsilon}$ agrees
    with the scaled Lawlor potential when
    $|\mathbf{x}|\leq\varepsilon^{\zeta}$ and vanishes when
    $|\mathbf{x}|\geq 2\varepsilon^{\zeta}$. The inner part of
    $N^{\varepsilon}$ is taken directly from
    $\varepsilon\overline{L}$, and $N^{\varepsilon}$ agrees with $X$
    outside the gluing region. It follows from
    \eqref{eq: outer Lawlor potential asymptotic} that, for every
    integer $k\geq 0$, there is $C_k>0$, independent of
    $\varepsilon$, such that
    \begin{equation}\label{eq: cutoff potential decay}
        |\nabla^k\mathfrak{D}_{\varepsilon}(\mathbf{x})|
        \leq C_k\varepsilon^n|\mathbf{x}|^{2-n-k}
        \qquad\text{for }R_L\varepsilon<|\mathbf{x}|<r_0.
    \end{equation}

    \item[(b)] \emph{The neck scale.}
    For some sufficiently large $t_0$, the function
    $\varepsilon:[t_0,\infty)\to(0,\infty)$ is $C^1$ and satisfies
    \[
        D^{-1}t^{-\frac{1}{n-2}}
        \leq\varepsilon(t)\leq Dt^{-\frac{1}{n-2}},
        \qquad
        |\varepsilon'(t)|\leq Dt^{\frac{1-n}{n-2}}
    \]
    for some $D>1$.

    \item[(c)] \emph{The Hamiltonian perturbation.}
    The flow $N(t)$ is a Hamiltonian deformation of
    $N^{\varepsilon(t)}$ determined by a potential $u(\cdot,t)$
    satisfying
    \begin{equation}\label{eq: perturbation weighted bound}
        \sup_{t\geq t_0}\;
        \sup_{\mathbf{x}\in N^{\varepsilon(t)}}
        \sum_{j=0}^{2}
        t^{\mu}\rho_{\varepsilon(t)}(\mathbf{x})^{\nu+j}
        |\nabla^ju(\mathbf{x},t)|\leq 1,
    \end{equation}
    where
    \[
        0<\nu<n-2,
        \qquad \mu>\frac{\nu+2}{n-2}.
    \]
    Here the derivatives are taken with respect to the induced metric
    on $N^{\varepsilon(t)}$, and the weight satisfies
    \[
        \rho_{\varepsilon(t)}(\mathbf{x})
        \sim\sqrt{\varepsilon(t)^2+|\mathbf{x}|^2}
        \quad\text{in }B_{r_0}(0),
        \qquad
        \rho_{\varepsilon(t)}\sim 2r_0
        \quad\text{outside }B_{2r_0}(0).
    \]
\end{enumerate}

On the graphical region, we also write $u$ for its expression in the
cone coordinates of the gluing construction. In these coordinates,
$N(t)$ has Lagrangian graph potential
$\mathfrak{D}_{\varepsilon(t)}+u(\cdot,t)$.
On annuli where the approximate solution is uniformly graphical, the
weighted estimates above also hold for the coordinate derivatives, up
to uniform constants. We use the graphical parametrisations to
transport the weight to $N(t)$ as well.

\begin{proposition}
    Let $\{N(t)\}_{t\geq t_0}$ be the flow in
    Theorem~\ref{thm: STW24 solution}. Set
    \[
        R(t):=r_0,
        \qquad \lambda(t):=\varepsilon(t)
        \sim t^{-\frac{1}{n-2}}.
    \]
    Then $(N(t),R(t),\lambda(t))$ satisfies
    Assumption~\ref{blow-up assumptions} in the flat coordinates
    centred at $x_*$.
\end{proposition}

\begin{proof}
    The outer scale $R$ is constant, hence $C^1$ and non-increasing,
    and the inner scale $\lambda$ is $C^1$ by (b). Moreover,
    $\lambda(t)/R(t)\to 0$, which gives
    \eqref{scale assumption 1}, while
    \[
        R(t)^2\frac{|\lambda'(t)|}{\lambda(t)}
        \lesssim r_0^2t^{-1}\longrightarrow 0
    \]
    gives \eqref{scale assumption 2}.
    The convergence at the neck scale in
    \eqref{assumption, blow-up to Lawlor neck} follows from
    Theorem~\ref{thm: STW24 solution} and the scale estimates in (b).

    It remains to verify
    \eqref{assumption, blow-up to cone}.
    Let $t_i\to\infty$ and $\mu_i\to 0$ satisfy
    \[
        \varepsilon(t_i)\ll\mu_i\ll r_0,
    \]
    and consider the parabolically rescaled flows
    \[
        N_i(\tau)
        :=\mu_i^{-1}\bigl(N(t_i+\mu_i^2\tau)
        \cap B_{2r_0}(0)\bigr).
    \]
    Write $s_i(\tau):=t_i+\mu_i^2\tau$ and
    $\varepsilon_i(\tau):=\varepsilon(s_i(\tau))$.
    Fix $\Lambda>1$ and restrict to $|\tau|\leq\Lambda$.
    Since $|\varepsilon'(t)|/\varepsilon(t)\lesssim t^{-1}$,
    \[
        \sup_{|\tau|\leq\Lambda}
        \left|\log\frac{\varepsilon_i(\tau)}{\varepsilon(t_i)}\right|
        \lesssim_{\Lambda}\frac{\mu_i^2}{t_i}
        \longrightarrow 0.
    \]
    In particular,
    \begin{equation}\label{eq: inner scale comparison}
        \frac{\varepsilon_i(\tau)}{\mu_i}
        =\frac{\varepsilon_i(\tau)}{\varepsilon(t_i)}
        \frac{\varepsilon(t_i)}{\mu_i}
        \longrightarrow 0
    \end{equation}
    uniformly for $|\tau|\leq\Lambda$.

    For all sufficiently large $i$, the rescaled flow is therefore
    graphical over the fixed annulus
    \[
        A_{\Lambda}
        :=\mathcal{C}\cap
        \{\Lambda^{-1}\leq|\mathbf{x}|\leq\Lambda\},
    \]
    with potential
    \[
        f_i(\mathbf{x},\tau)
        :=\mu_i^{-2}\bigl[
        \mathfrak{D}_{\varepsilon_i(\tau)}(\mu_i\mathbf{x})
        +u(\mu_i\mathbf{x},s_i(\tau))\bigr].
    \]
    For $k=0,1,2$, \eqref{eq: cutoff potential decay} gives
    \begin{align*}
        \left|\nabla^k\left[
        \mu_i^{-2}\mathfrak{D}_{\varepsilon_i(\tau)}
        (\mu_i\mathbf{x})\right]\right|
        &\leq C_{\Lambda}
        \mu_i^{k-2}\varepsilon_i(\tau)^n\mu_i^{2-n-k}\\
        &=C_{\Lambda}
        \left(\frac{\varepsilon_i(\tau)}{\mu_i}\right)^n
        \longrightarrow 0.
    \end{align*}
    On the same annulus,
    $\rho_{\varepsilon_i(\tau)}(\mu_i\mathbf{x})\sim\mu_i$
    uniformly in $i$ and $\tau$. Thus
    \eqref{eq: perturbation weighted bound} yields
    \begin{align*}
        \left|\nabla^k\left[
        \mu_i^{-2}u(\mu_i\mathbf{x},s_i(\tau))\right]\right|
        &\leq C_{\Lambda}
        s_i(\tau)^{-\mu}\mu_i^{-\nu-2}\\
        &=C_{\Lambda}s_i(\tau)^{-\mu}
        \varepsilon_i(\tau)^{-\nu-2}
        \left(\frac{\varepsilon_i(\tau)}{\mu_i}\right)^{\nu+2}\\
        &\lesssim_{\Lambda}
        s_i(\tau)^{-\mu+\frac{\nu+2}{n-2}}
        \longrightarrow 0,
    \end{align*}
    where the last step uses
    $\mu>(\nu+2)/(n-2)$.
    All these limits are uniform on
    $A_{\Lambda}\times[-\Lambda,\Lambda]$.

    Consequently, $N_i(\tau)$ converges to the stationary cone
    $\mathcal{C}$ locally in the spatial $C^1$ topology, uniformly on
    compact time intervals away from the origin.
    On each compact annulus, the graphical mean curvature flow
    equation is uniformly parabolic for all sufficiently large $i$.
    White's local regularity theorem \cite{White05} and interior estimates for mean curvature flow then upgrade
    this convergence to smooth convergence on smaller parabolic
    cylinders. Since $\Lambda$ is arbitrary, this proves
    \eqref{assumption, blow-up to cone}.
\end{proof}

\begin{proposition}\label{prop: weighted angle bound in annulus infinite time}
    Let $\{N(t)\}_{t\geq t_0}$ be the flow in
    Theorem~\ref{thm: STW24 solution}, and normalise its Lagrangian
    angle $\theta$ to converge to $0$ away from $x_*$.
    For every $\alpha\in\mathbb{R}$,
    $\beta\leq\mu(n-2)$, and $C>1$, there exist $t_*\geq t_0$ and
    $M>0$ such that
    \[
        \sup_{t\geq t_*}\;
        \sup_{\substack{\mathbf{x}\in N(t)\\
        C^{-1}r_0\leq|\mathbf{x}|\leq r_0}}
        \frac{\rho_{\varepsilon(t)}(\mathbf{x})^{\alpha}}
        {\varepsilon(t)^{\beta}}
        |\theta(\mathbf{x},t)|\leq M.
    \]
\end{proposition}

\begin{proof}
    Since $\varepsilon(t)\to 0$, the gluing region is eventually
    contained strictly inside $B_{C^{-1}r_0}(0)$.
    By increasing $t_*$ if necessary, the portion of $N(t)$ in
    $\{C^{-1}r_0\leq|\mathbf{x}|\leq r_0\}$ is therefore the
    Lagrangian graph of $\nabla u(\cdot,t)$ over the two planes of
    $\mathcal{C}$ for all $t\geq t_*$.
    In these graphical coordinates,
    \[
        \theta(\cdot,t)
        =\sum_{j=1}^{n}
        \arctan\lambda_j\bigl(\nabla^2u(\cdot,t)\bigr),
    \]
    where $\lambda_j(\nabla^2u)$ are the eigenvalues of
    $\nabla^2u$. Hence
    \[
        |\theta(\cdot,t)|
        \leq\sum_{j=1}^{n}
        |\lambda_j(\nabla^2u(\cdot,t))|
        \leq\sqrt{n}\,|\nabla^2u(\cdot,t)|.
    \]
    Combining this with
    \eqref{eq: perturbation weighted bound}, we obtain
    \[
        \frac{\rho_{\varepsilon(t)}^{\alpha}}
        {\varepsilon(t)^{\beta}}|\theta|
        \lesssim
        t^{-\mu}\varepsilon(t)^{-\beta}
        \rho_{\varepsilon(t)}^{\alpha-\nu-2}.
    \]
    On this fixed annulus, the weight is bounded above and below by
    positive constants independent of $t$, so
    $\rho_{\varepsilon(t)}^{\alpha-\nu-2}$ is uniformly bounded
    for every $\alpha\in\mathbb{R}$. Moreover, (b) implies
    \[
        t^{-\mu}\varepsilon(t)^{-\beta}
        \lesssim t^{-\mu+\frac{\beta}{n-2}}.
    \]
    The exponent is non-positive because $\beta\leq\mu(n-2)$,
    which proves the claimed bound.
\end{proof}

\begin{corollary} \label{main cor infinite-time sings sect}
    Let $\{N(t)\}_{t\geq t_0}$ be the flow in
    Theorem~\ref{thm: STW24 solution}. Then
    \[
        \sup_{t\geq t_0}\;\sup_{\mathbf{x}\in N(t)}
        |\mathbf{H}_{N(t)}(\mathbf{x})|<\infty.
    \]
\end{corollary}

\begin{proof}
    Since $\mu(n-2)>\nu+2>2$, we may choose
    \[
        0<\alpha<\min\{\mu(n-2)-1,n-2\},
        \qquad \beta:=1+\alpha<\mu(n-2).
    \]
    Let $C>1$ be the constant in
    Corollary~\ref{cor: pointwise mean curvature bound}.
    The preceding propositions and that corollary give constants
    $M>0$ and $t_*\geq t_0$ such that
    \[
        |\mathbf{H}_{N(t)}(\mathbf{x})|
        \leq C(1+M)
        \qquad\text{for }|\mathbf{x}|\leq C^{-1}r_0,
        \quad t\geq t_*.
    \]
    On the complement of this ball,
    Theorem~\ref{thm: STW24 solution} gives smooth convergence to
    $X_1\cup X_2$, so the mean curvature converges uniformly to $0$
    there. Combining these bounds proves the result.
\end{proof}

\bibliography{BibFile}

\end{document}